\documentclass[A4j,11pt]{article}
\usepackage{graphics}
\usepackage{amsmath}
\usepackage{amssymb}
\usepackage{latexsym}
\usepackage{amsfonts}
\begin{document}

\title{{\bf Examples of a complex hyperpolar action\\
without singular orbit}}
\author{{\bf Naoyuki Koike}}
%
%
\date{}
%
\maketitle
\begin{abstract}
The notion of a complex hyperpolar action on a symmetric space of non-compact 
type has recently been introduced as a counterpart to the hyperpolar action on 
a symmetric space of compact type.  As examples of a complex hyperpolar 
action, we have Hermann type actions, which admit a totally geodesic singular 
orbit (or a fixed point) except for one example.  All principal orbits of 
Hermann type actions are curvature-adapted and proper complex equifocal.  
In this paper, we give some examples of a complex hyperpolar action without 
singular orbit as solvable group free actions and find complex hyperpolar 
actions all of whose orbits are non-curvature-adapted or non-proper complex 
equifocal among the examples.  Also, we show that some of the examples 
possess the only minimal orbit.  
\end{abstract}

\vspace{0.5truecm}

\section{Introduction}
In symmetric spaces, the notion of an equifocal submanifold was introduced 
in [TT].  This notion is defined as a compact submanifold with globally flat 
and abelian normal bundle such that the focal radius functions for each 
parallel normal vector field are constant.  
However, this conditions of the equifocality is rather weak in the case 
where the symmetric spaces are of non-compact type and the submanifold is 
non-compact.  So we [Koi1,2] 
have recently introduced the notion of a complex equifocal 
submanifold in a symmetric space $G/K$ of non-compact type.  
This notion is defined by imposing the constancy of the complex focal 
radius functions instead of focal radius functions.  
Here we note that the complex focal radii are the 
quantities indicating the positions of the focal points of the 
extrinsic complexification of the submanifold, where the submanifold needs 
to be assumed to be complete and of class $C^{\omega}$ 
(i.e., real analytic).  
On the other hand, Heintze-Liu-Olmos [HLO] has recently defined the 
notion of an isoparametric submanifold with flat section in a general 
Riemannian manifold as a submanifold such that the normal holonomy group is 
trivial, its sufficiently close parallel submanifolds are of constant mean 
curvature with respect to the radial direction and that the image of the 
normal space at each point by the normal 
exponential map is flat and totally geodesic.  
We [Koi2] showed the following fact:

\vspace{0.3truecm}

{\sl All isoparametric submanifolds with 
flat section in a symmetric space $G/K$ of non-compact type are complex 
equifocal and that conversely, all curvature-adapted and complex equifocal 
submanifolds are isoparametric ones with flat section.}

\vspace{0.3truecm}

Here the curvature-adaptedness means that, for each normal vector $v$ of 
the submanifold, the Jacobi operator $R(\cdot,v)v$ preserves the tangent space 
of the submanifold invariantly and the restriction of $R(\cdot,v)v$ to the 
tangent space commutes with the shape operator $A_v$, where $R$ is the 
curvature tensor of $G/K$.  
Furthermore, as a subclass of the class of complex equifocal submanifolds, 
we [Koi3] defined that of the proper complex equifocal submanifolds in $G/K$ 
as a complex equifocal submanifold whose lifted submanifold to 
$H^0([0,1],\mathfrak g)$ ($\mathfrak g:={\rm Lie}\,G$) through 
some pseudo-Riemannian submersion of $H^0([0,1],\mathfrak g)$ onto $G/K$ is 
proper complex isoparametric in the sense of [Koi1], where we note that 
$H^0([0,1],\mathfrak g)$ is a pseudo-Hilbert space consisting of certain kind 
of paths in the Lie algebra $\mathfrak g$ of $G$.  
Let $G/K$ be a symmetric space of non-compact type and $H$ be a closed 
subgroup of $G$ which admits an embedded complete flat submanifold meeting all 
$H$-orbits orthogonally.  Then the $H$-action on $G/K$ is called a 
{\it complex hyperpolar action}.  
This action was named thus because this action has not necessarily a singular 
orbit (which should be called a pole of this action) but the complexified 
action has a singular orbit.  Note that all cohomogeneity one actions are 
complex hyperpolar.  
We [Koi2] showed that principal orbits of 
a complex hyperpolar actions are isoparametric submanifolds with flat section 
and hence they are complex equifocal.  Conversely we 
[Koi5] have recently showed that all homogeneous complex equifocal 
submanifolds occurs as principal orbits of complex hyperpolar actions.  Let 
$H'$ be a symmetric subgroup of $G$ (i.e., there exists an involution $\sigma$ 
of $G$ with $({\rm Fix}\,\sigma)_0\subset H'\subset{\rm Fix}\,\sigma$), where 
${\rm Fix}\,\sigma$ is the fixed point group of $\sigma$ and 
$({\rm Fix}\,\sigma)_0$ is the identity component of ${\rm Fix}\,\sigma$.  
Then the $H'$-action on $G/K$ is called a {\it Hermann type action}.  
A Hermann type action admits a totally geodesic orbit or a fixed point.  
Except for one example, the totally geodsic orbit is singular (see Theorem E 
of [Koi5]).  We [Koi3] showed that principal orbits of 
a Hermann type action are proper complex equifocal and curvature-adapted.  
We [Koi5] have recently showed that all complex hyperpolar actions of 
cohomogeneity greater than one on $G/K$ admitting a totally geodesic orbit and 
all complex hyperpolar actions of cohomogeneity one on $G/K$ admitting 
reflective orbit are orbit equivalent to Hermann type actions (see Theorems B, 
C and Remark 1.1 in [Koi5]).  
Let $G/K$ be a symmetric space of non-compact type, $\mathfrak g=\mathfrak f
+\mathfrak p$ ($\mathfrak f:={\rm Lie}\,K$) be the Cartan decomposition 
associated with $(G,K)$, $\mathfrak a$ be the maximal abelian subspace of 
$\mathfrak p$, $\widetilde{\mathfrak a}$ be the Cartan subalgebra of 
$\mathfrak g$ containing $\mathfrak a$ and $\mathfrak g=\mathfrak f+\mathfrak a
+\mathfrak n$ be the Iwasawa's decomposition.  Let $A,\,\widetilde A$ and $N$ 
be the connected Lie subgroups of $G$ having $\mathfrak a,\,
\widetilde{\mathfrak a}$ and $\mathfrak n$ as their Lie algebras, 
respectively.  
Let $\pi:G\to G/K$ be the natural projection.  
The symmetric space $G/K$ is 
identified with the solvable group $AN$ with a left-invariant metric 
through $\pi\vert_{AN}$.  
In this paper, we first prove the following fact for a complex hyperpolar 
action without singular orbit.  

\vspace{0.5truecm}

\noindent
{\bf Theorem A.} {\sl Any complex hyperpolar action on $G/K(=AN)$ without 
singular orbit is orbit equivalent to the free action of some solvable group 
contained in $\widetilde AN$.}

\vspace{0.5truecm}

Next we give some examples of a complex hyperpolar 
action without singular orbit as the free actions of solvable groups contained 
in $AN$ (see Examples 1 and 2 of Section 3), which contain examples of 
cohomogeneity one actions without singular orbit constructed by J. Berndt and 
H. Tamaru [BT1] as special cases (see also [B]).  
Among these examples, we find complex hyperpolar actions all of whose orbits 
are non-proper complex equifocal or non-curvature-adapted.  
As its result, we have the following facts.  

\vspace{0.5truecm}

\noindent
{\bf Theorem B.} {\sl {\rm (i)} For any symmetric space $G/K$ of 
non-compact type and any positive integer $r$ with $r\leq {\rm rank}(G/K)$, 
there exists a complex hyperpolar action without singular orbit 
such that the cohomogeneity is equal to $r$ and that any of the orbits is 
not proper complex equifocal.  

{\rm (ii)} Let $G/K$ be one of $SU(p,q)/S(U(p)\times U(q))\,\,(p<q),\,\,
Sp(p,q)/Sp(p)\times Sp(q)\,\,(p<q),\,\,SO^{\ast}(2n)/U(n)\,\,(n:{\rm odd}),
\,\,E_6^{-14}/Spin(10)\cdot U(1)$ or $F_4^{-20}/Spin(9)$.  
Then, for any positive integer $r$ with $r\leq {\rm rank}(G/K)$, there exists 
a complex hyperpolar action without singular orbit 
such that the cohomogeneity is equal to $r$ and that any of the orbits is 
not curvature-adapted.  
}

\vspace{0.5truecm}

Also, among those examples, we find complex hyperpolar actions possessing the 
only minimal orbit.  As its result, we have the following fact.  

\vspace{0.5truecm}

\noindent
{\bf Theorem C.} {\sl For any irreducible symmetric space $G/K$ of 
non-compact type and any positive integer $r\leq[\frac12({\rm rank}(G/K)+1)]$, 
there exists a complex hyperpolar action without singular orbit 
such that the cohomogeneity is equal to $r$ and that the only orbit is 
minimal.}

\section{Complex equifocal submanifolds}
In this section, we recall the notions of a complex equifocal submanifold 
and a proper complex equifocal submanifold.  
We first recall the notion of a complex equifocal submanifold.  
Let $M$ be an immersed submanifold with abelian normal bundle 
in a symmetric space $N=G/K$ of non-compact type.  Denote by $A$ the shape 
tensor of $M$.  Let $v\in T^{\perp}_xM$ and $X\in T_xM$ ($x=gK$).  Denote 
by $\gamma_v$ the geodesic in $N$ with $\dot{\gamma}_v(0)=v$.  
The strongly $M$-Jacobi field $Y$ along $\gamma_v$ with $Y(0)=X$ (hence 
$Y'(0)=-A_vX$) is given by 
$$Y(s)=(P_{\gamma_v\vert_{[0,s]}}\circ(D^{co}_{sv}-sD^{si}_{sv}\circ A_v))
(X),$$
where $Y'(0)=\widetilde{\nabla}_vY,\,\,P_{\gamma_v\vert_{[0,s]}}$ is 
the parallel translation along $\gamma_v\vert_{[0,s]}$ and 
$D^{co}_{sv}$ (resp. $D^{si}_{sv}$) is given by 
$$\begin{array}{c}
\displaystyle{
D^{co}_{sv}=g_{\ast}\circ\cos(\sqrt{-1}{\rm ad}(sg_{\ast}^{-1}v))
\circ g_{\ast}^{-1}}\\
\displaystyle{\left({\rm resp.}\,\,\,\,
D^{si}_{sv}=g_{\ast}\circ
\frac{\sin(\sqrt{-1}{\rm ad}(sg_{\ast}^{-1}v))}
{\sqrt{-1}{\rm ad}(sg_{\ast}^{-1}v)}\circ g_{\ast}^{-1}\right).}
\end{array}$$ 
Here ${\rm ad}$ is the adjoint 
representation of the Lie algebra $\mathfrak g$ of $G$.  
All focal radii of $M$ along $\gamma_v$ are obtained as real numbers $s_0$ 
with ${\rm Ker}(D^{co}_{s_0v}-s_0D^{si}_{s_0v}\circ A_v)\not=\{0\}$.  So, we 
call a complex number $z_0$ with ${\rm Ker}(D^{co}_{z_0v}-
z_0D^{si}_{z_0v}\circ A_v^{{\bf c}})\not=\{0\}$ a {\it complex 
focal radius of} $M$ {\it along} $\gamma_v$ and call ${\rm dim}\,
{\rm Ker}(D^{co}_{z_0v}-z_0D^{si}_{z_0v}\circ A_v^{{\bf c}})$ the 
{\it multiplicity} of the complex focal radius $z_0$, 
where $A_v^{\bf c}$ is the complexification of $A_v$ and $D^{co}_{z_0v}$ 
(resp. $D^{si}_{z_0v}$) is a ${\bf C}$-linear transformation of 
$(T_xN)^{\bf c}$ defined by 
$$\begin{array}{c}
\displaystyle{
D^{co}_{z_0v}=g^{\bf c}_{\ast}\circ\cos(\sqrt{-1}{\rm ad}^{\bf c}
(z_0g_{\ast}^{-1}v))\circ (g^{\bf c}_{\ast})^{-1}}\\
\displaystyle{\left({\rm resp.}\,\,\,\,
D^{si}_{sv}=g^{\bf c}_{\ast}\circ
\frac{\sin(\sqrt{-1}{\rm ad}^{\bf c}(z_0g_{\ast}^{-1}v))}
{\sqrt{-1}{\rm ad}^{\bf c}(z_0g_{\ast}^{-1}v)}\circ(g^{\bf c}_{\ast})^{-1}
\right),}
\end{array}$$
where $g_{\ast}^{\bf c}$ (resp. ${\rm ad}^{\bf c}$) is the complexification 
of $g_{\ast}$ (resp. ${\rm ad}$).  
Here we note that, in the case where $M$ is of class $C^{\omega}$, 
complex focal radii along $\gamma_v$ 
indicate the positions of focal points of the extrinsic 
complexification $M^{\bf c}(\hookrightarrow G^{\bf c}/K^{\bf c})$ of $M$ 
along the complexified geodesic $\gamma_{\iota_{\ast}v}^{\bf c}$, where 
$G^{\bf c}/K^{\bf c}$ is the anti-Kaehlerian symmetric space associated with 
$G/K$ and $\iota$ is the natural immersion of $G/K$ into 
$G^{\bf c}/K^{\bf c}$.  
See Section 4 of [Koi2] about the definitions of 
$G^{\bf c}/K^{\bf c},\,
M^{\bf c}(\hookrightarrow G^{\bf c}/K^{\bf c})$ and 
$\gamma_{\iota_{\ast}v}^{\bf c}$.  
Also, for a complex focal radius $z_0$ of $M$ along $\gamma_v$, we 
call $z_0v$ ($\in (T^{\perp}_xM)^{\bf c}$) a 
{\it complex focal normal vector of} $M$ {\it at} $x$.  
Furthermore, assume that $M$ has globally flat normal bundle, that is, 
the normal holonomy group of $M$ is trivial.  
Let $\tilde v$ be a parallel unit normal vector field of $M$.  
Assume that the number (which may be $0$ and $\infty$) of distinct complex 
focal radii along $\gamma_{\tilde v_x}$ is independent of the choice of 
$x\in M$.  Furthermore assume that the number is not equal to $0$.  
Let $\{r_{i,x}\,\vert\,i=1,2,\cdots\}$ 
be the set of all complex focal radii along $\gamma_{\tilde v_x}$, where 
$\vert r_{i,x}\vert\,<\,\vert r_{i+1,x}\vert$ or 
"$\vert r_{i,x}\vert=\vert r_{i+1,x}\vert\,\,\&\,\,{\rm Re}\,r_{i,x}
>{\rm Re}\,r_{i+1,x}$" or 
"$\vert r_{i,x}\vert=\vert r_{i+1,x}\vert\,\,\&\,\,
{\rm Re}\,r_{i,x}={\rm Re}\,r_{i+1,x}\,\,\&\,\,
{\rm Im}\,r_{i,x}=-{\rm Im}\,r_{i+1,x}<0$".  
Let $r_i$ ($i=1,2,\cdots$) be complex 
valued functions on $M$ defined by assigning $r_{i,x}$ to each $x\in M$.  
We call these functions $r_i$ ($i=1,2,\cdots$) {\it complex 
focal radius functions for} $\tilde v$.  
We call $r_i\tilde v$ a {\it complex focal normal vector field for} 
$\tilde v$.  If, for each parallel 
unit normal vector field $\tilde v$ of $M$, the number of distinct complex 
focal radii along $\gamma_{\tilde v_x}$ is independent of the choice of 
$x\in M$, each complex focal radius function for $\tilde v$ 
is constant on $M$ and it has constant multiplicity, then 
we call $M$ a {\it complex equifocal submanifold}.  

Let $N=G/K$ be a symmetric space of non-compact type and $\pi$ be the 
natural projection of $G$ onto $G/K$.  
Let $(\mathfrak g,\theta)$ be the orthogonal symmetric Lie algebra of $G/K$, 
$\mathfrak f=\{X\in \mathfrak g\,\vert\,\theta(X)=X\}$ and 
$\mathfrak p=\{X\in \mathfrak g\,\vert\,\theta(X)=-X\}$, which is identified 
with the tangent space $T_{eK}N$.  Let $\langle\,\,,\,\,\rangle$ be 
the ${\rm Ad}(G)$-invariant non-degenerate symmetric bilinear form of 
$\mathfrak g$ inducing the Riemannian metric of $N$.  Note that 
$\langle\,\,,\,\,\rangle\vert_{\mathfrak f\times\mathfrak f}$ 
(resp. $\langle\,\,,\,\,\rangle\vert_{\mathfrak p\times\mathfrak p}$) is 
negative (resp. positive) definite.  Denote by the same symbol 
$\langle\,\,,\,\,\rangle$ the bi-invariant pseudo-Riemannian metric of $G$ 
induced from $\langle\,\,,\,\,\rangle$ and the Riemannian metric of $N$.  
Set $\mathfrak g_+:=\mathfrak p,\,\,\mathfrak g_-:=\mathfrak f$ and 
$\langle\,\,,\,\,\rangle_{\mathfrak g_{\pm}}
:=-\pi_{\mathfrak g_-}^{\ast}\langle\,\,,\,\,\rangle
+\pi_{\mathfrak g_+}^{\ast}\langle\,\,,\,\,\rangle$, where 
$\pi_{\mathfrak g_-}$ (resp. $\pi_{\mathfrak g_+}$) is the projection of 
$\mathfrak g$ onto $\mathfrak g_-$ (resp. $\mathfrak g_+$).  Let 
$H^0([0,1],\mathfrak g)$ be the space of all $L^2$-integrable paths 
$u:[0,1]\to\mathfrak g$ (with respect to 
$\langle\,\,,\,\,\rangle_{\mathfrak g_{\pm}}$).  
Let $H^0([0,1],\mathfrak g_-)$ (resp. $H^0([0,1],\mathfrak g_+)$) be 
the space of all $L^2$-integrable paths 
$u:[0,1]\to\mathfrak g_-$ (resp. $u:[0,1]\to \mathfrak g_+$) with respect to 
$-\langle\,\,,\,\,\rangle\vert_{\mathfrak g_-\times\mathfrak g_-}$ 
(resp. $\langle\,\,,\,\,\rangle\vert_{\mathfrak g_+\times\mathfrak g_+}$).  
It is clear that $H^0([0,1],\mathfrak g)=H^0([0,1],\mathfrak g_-)\oplus 
H^0([0,1],\mathfrak g_+)$.  
Define a non-degenerate symmetric bilinear form 
$\langle\,\,,\,\,\rangle_0$ of 
$H^0([0,1],\mathfrak g)$ by 
$\langle u,v\rangle_0:=\int_0^1\langle u(t),v(t)\rangle dt$.  It is easy to 
show that the decomposition $H^0([0,1],\mathfrak g)=H^0([0,1],\mathfrak g_-)
\oplus H^0([0,1],\mathfrak g_+)$ is an orthogonal time-space decomposition 
with respect to $\langle\,\,,\,\,\rangle_0$.  For simplicity, set 
$H^0_{\pm}:=H^0([0,1],\mathfrak g_{\pm})$ and 
$\langle\,\,,\,\,\rangle_{0,H^0_{\pm}}:=-\pi^{\ast}_{H^0_-}
\langle\,\,,\,\,\rangle_0+\pi^{\ast}_{H^0_+}\langle\,\,,\,\,\rangle_0$, where 
$\pi_{H^0_-}$ (resp. $\pi_{H^0_+}$) is the projection of 
$H^0([0,1],\mathfrak g)$ onto $H^0_-$ (resp. $H^0_+$).  
It is clear that $\langle u,v\rangle_{0,H^0_{\pm}}
=\int_0^1\langle u(t),v(t)\rangle_{\mathfrak g_{\pm}}dt$ 
($u,\,v\in H^0([0,1],\mathfrak g)$).  
Hence $(H^0([0,1],\mathfrak g),\,\langle\,\,,\,\,\rangle_{0,H^0_{\pm}})$ is a 
Hilbert space, that is, $(H^0([0,1],\mathfrak g),\,\langle\,\,,\,\,\rangle_0)$ 
is a pseudo-Hilbert space.  Let $H^1([0,1],G)$ be the Hilbert Lie group of all 
absolutely continuous paths $g:[0,1]\to G$ such that the weak derivative $g'$ 
of $g$ is squared integrable (with respect to 
$\langle\,\,,\,\,\rangle_{\mathfrak g_{\pm}}$), that is, 
$g_{\ast}^{-1}g'\in H^0([0,1],\mathfrak g)$.  Define a map 
$\phi:H^0([0,1],\mathfrak g)\to G$ by $\phi(u)=g_u(1)$ 
($u\in H^0([0,1],\mathfrak g)$), where $g_u$ is the element of 
$H^1([0,1],G)$ satisfying $g_u(0)=e$ and $g_{u\ast}^{-1}g_u'=u$.  
We call this map the {\it parallel transport map} (from $0$ to $1$).  
This submersion $\phi$ is a pseudo-Riemannian submersion of 
$(H^0([0,1],\mathfrak g),\langle\,\,,\,\,\rangle_0)$ onto 
$(G,\langle\,\,,\,\,\rangle)$.  
Let $\pi:G\to G/K$ be the natural projection.  It follows from Theorem A of 
[Koi1] (resp. Theorem 1 of [Koi2]) that, in the case where $M$ is curvature 
adapted (resp. of class $C^{\omega}$), $M$ is complex equifocal if and only if 
each component of $(\pi\circ\phi)^{-1}(M)$ is complex isoparametric.  
See [Koi1] about the definition of a complex isoparametric submanifold in a 
pseudo-Hilbert space.  
In particular, if each component of $(\pi\circ\phi)^{-1}(M)$ 
are proper complex isoparametric, that is, for each normal vector $v$, 
there exists a pseudo-orthonormal base of the complexified tangent sapce 
consisting of the eigenvectors of the complexified shape operator for $v$, 
then we call $M$ a {\it proper complex equifocal submanifold}.  

Next we recall the notion of an infinite dimensional 
anti-Kaehlerian isoparametric submanifold.  
Let $M$ be an anti-Kaehlerian Fredholm submanifold in an infinite 
dimensional anti-Kaehlerian space $V$ and $A$ be the shape tensor of $M$.  
See [Koi2] about the definitions of an infinite dimensional 
anti-Kaehlerian space and anti-Kaehlerian Fredholm submanifold in the space.  
Denote by the same symbol $J$ the complex structures of $M$ and $V$.  
Fix a unit normal vector $v$ of $M$.  If there exists $X(\not=0)\in TM$ with 
$A_vX=aX+bJX$, then we call the complex number $a+b\sqrt{-1}$ a 
$J$-{\it eigenvalue of} $A_v$ (or a {\it complex principal curvature of 
direction} $v$) and call $X$ a $J$-{\it eigenvector for} $a+b\sqrt{-1}$.  
Also, we call the space of all $J$-eigenvectors for 
$a+b\sqrt{-1}$ a $J$-{\it eigenspace for} $a+b\sqrt{-1}$.  
The $J$-eigenspaces are orthogonal to one another and 
each $J$-eigenspace is $J$-invariant.  
We call the set of all $J$-eigenvalues of $A_v$ the $J$-{\it spectrum of} 
$A_v$ and denote it by ${{\rm Spec}}_JA_v$.  
The set ${{\rm Spec}}_JA_v\setminus\{0\}$ 
is described as follows:
$${{\rm Spec}}_JA_v\setminus\{0\}
=\{\lambda_i\,\vert\,i=1,2,\cdots\}$$
$$\left(
\begin{array}{c}
\displaystyle{\vert\lambda_i\vert>\vert\lambda_{i+1}\vert\,\,\,{{\rm or}}
\,\,\,{\rm "}\vert\lambda_i\vert=\vert\lambda_{i+1}\vert\,\,\&\,\,
{{\rm Re}}\,\lambda_i>{{\rm Re}}\,\lambda_{i+1}{\rm "}}\\
\displaystyle{{{\rm or}}\,\,\,"\vert\lambda_i\vert=\vert\lambda_{i+1}\vert 
\,\,\&\,\,
{{\rm Re}}\,\lambda_i={{\rm Re}}\,\lambda_{i+1}\,\,\&\,\,
{{\rm Im}}\,\lambda_i=-{{\rm Im}}\,\lambda_{i+1}>0"}
\end{array}
\right).$$
Also, the $J$-eigenspace for each $J$-eigenvalue of $A_v$ other than $0$ 
is of finite dimension.  
We call the $J$-eigenvalue $\lambda_i$ 
the $i$-{\it th complex principal curvature of direction} $v$.  Assume that 
$M$ has globally flat normal bundle.  Fix a 
parallel normal vector field $\widetilde v$ of $M$.  
Assume that the number 
(which may be $\infty$) of distinct complex principal curvatures of direction 
$\widetilde v_x$ is independent of the choice of $x\in M$.  
Then we can define functions $\widetilde{\lambda}_i$ ($i=1,2,\cdots$) on $M$ 
by assigning the $i$-th complex principal curvature of direction 
$\widetilde v_x$ to each $x\in M$.  We call this function 
$\widetilde{\lambda}_i$ the $i$-{\it th complex principal curvature 
function of direction} $\widetilde v$.  
If $M$ satisfies the following condition (AKI), then we call 
$M$ an {\it anti-Kaehlerian isoparametric submanifold}:

\vspace{0.2truecm}

\noindent
(AKI) For each parallel normal vector field $\widetilde v$, the 
number of distinct complex principal curvatures of direction $\widetilde v_x$ 
is independent of the choice of $x\in M$, each complex principal curvature 
function of direction $\widetilde v$ is constant on $M$ and it has constant 
multiplicity.  

\vspace{0.2truecm}

\noindent
Let $\{e_i\}_{i=1}^{\infty}$ be an orthonormal system of $T_xM$.  If 
$\{e_i\}_{i=1}^{\infty}\cup\{Je_i\}_{i=1}^{\infty}$ is an orthonormal base 
of $T_xM$, then we call $\{e_i\}_{i=1}^{\infty}$ a $J$-{\it orthonormal base}. 
If there exists a $J$-orthonormal base consisting of $J$-eigenvectors of 
$A_v$, then $A_v$ is said to {\it be diagonalized with respect to the} 
$J$-{\it orthonormal base}.  If $M$ is anti-Kaehlerian isoparametric and, 
for each $v\in T^{\perp}M$, the shape operator $A_v$ is 
diagonalized with respect to a $J$-orthonormal base, then we call 
$M$ a {\it proper anti-Kaehlerian isoparametric submanifold}.  
For arbitrary two unit normal vector $v_1$ and $v_2$ of a proper 
anti-Kaehlerian isoparametric submanifold, the shape operators $A_{v_1}$ 
and $A_{v_2}$ are simultaneously diagonalized with respect to a 
$J$-orthonormal base.  
Let $M$ be a proper anti-Kaehlerian isoparametric submanifold in 
an infinite dimensional anti-Kaehlerian space $V$.  
Let $\{E_i\,\vert\,i\in I\}$ 
be the family of distributions on $M$ such that, for each $x\in M$, 
$\{E_i(x)\,\vert\,i\in I\}$ is the set of all common $J$-eigenspaces of 
$A_v$'s ($v\in T^{\perp}_xM$).  The relation 
$T_xM=\overline{\displaystyle{\mathop{\oplus}_{i\in I}E_i}}$ holds.  
Let $\lambda_i$ ($i\in I$) be the section 
of $(T^{\perp}M)^{\ast}\otimes {{\bf C}}$ such that $A_v={{\rm Re}}
\lambda_i(v){{\rm id}}+{{\rm Im}}\lambda_i(v)J$ on $E_i(\pi(v))$ 
for each $v\in T^{\perp}M$, where $\pi$ is the bundle projection of 
$T^{\perp}M$.  We call $\lambda_i$ ($i\in I$) {\it complex principal 
curvatures} of $M$ and call distributions 
$E_i$ ($i\in I$) {\it complex curvature distributions} of $M$.  

In the case where $M$ is a real analytic submanifold in a symmetric space 
$G/K$ of non-compact type, it is shown that $M$ is complex equifocal if and 
only if $(\pi^{\bf c}\circ\phi^{\bf c})^{-1}(M^{\bf c})$ is 
anti-Kaehlerian isoparametric, where $\pi^{\bf c}$ is 
the natural projection of $G^{\bf c}$ onto $G^{\bf c}/K^{\bf c}$ and 
$\phi^{\bf c}$ is the parallel transport map for $G^{\bf c}$ 
(which is defined in similar to the above $\phi$).  
Also, it is shown that $M$ is proper complex equifocal if and only if 
$(\pi^{\bf c}\circ\phi^{\bf c})^{-1}(M^{\bf c})$ is proper anti-Kaehlerian 
isoparametric.  

\section{Proof of Theorems A and B}
In this section, we first prove Theorem A.  

\vspace{0.5truecm}

\noindent{\it Proof of Theorem A.} 
Let $H$ be a complex hyperpolar action on $G/K(=AN)$ without singular orbit, 
$H=LR$ ($L\,:\,$ semi-simple, $R\,:\,$ solvable) 
be the Levi decomposition of $H$ and $L=K_LA_LN_L$ ($K_L\,:\,$
compact, $A_L\,:\,$abelian, $N_L\,:\,$nilpotent) be the Iwasawa 
decomposition of $L$.  
Since $K_L$ is compact, it has a fixed point $p_0$ by the Cartan's fixed point 
theorem.  Suppose that $K_L\cdot p\not\subset A_LN_LR\cdot p$ for some 
$p\in G/K$.  Then we have ${\rm dim}\,H\cdot p_0\,<\,{\rm dim}\,H\cdot p$, 
which implies that $H\cdot p_0$ is a singular orbit.  
This contradicts the 
fact that the $H$-action has no singular orbit.  Hence it follows that 
$K_L\cdot p\subset A_LN_LR\cdot p$ for any $p\in G/K$.  Therefore 
we can show that the $A_LN_LR$-action has the same orbits as the $H$-action.  
The group $A_LN_LR$ is decomposed into the product of some compact subgroup 
$T'$ and some solvable normal subgroup $S'$ admitting a maximal compact normal 
subgroup $S'_K$ contained in the center of $S'$ such that $S'/S'_K$ is simply 
connected (see Theorem 6 of [Ma]).  
Since $T'$ is compact, it is shown by the same argument as above that the 
$S'$-action has the same orbit as the $A_LN_LR$-action (hence the 
$H$-action).  Take any $p\in G/K$ and any $g\in S'$ with $g\not=e$.  Since 
$S'$ acts on $G/K$ effectively, there exists $p_1\in G/K$ with 
$g(p_1)\not=p_1$.  The section $\Sigma_{p_1}$ through $p_1$ is mapped into the 
section $\Sigma_{g(p_1)}$ through $g(p_1)$ by $g$.  
Since the $S'$-action has no singular 
orbit, we have $\Sigma_{p_1}\cap\Sigma_{g(p_1)}=\emptyset$.  Let $q$ be the 
intersection of $H\cdot p$ with $\Sigma_{p_1}$.  Then $g(q)$ is the 
intersection of $H\cdot p$ with $\Sigma_{g(p_1)}$.  Hence we have 
$g(q)\not= q$.  Therefore $S'$ acts on each $H$-orbit effectively.  Since the 
isotropy group $S'_p$ of $S'$ at any $p\in G/K$ is compact, it is contained 
in a conjugate of $S'_K$ (see Theorem 4 of [Ma]).  
Hence $S'_p$ is contained in the center of $S'$.  
Therefore, since the $S'_p$-action has a fixed point $p$ and it is effective, 
it is trivial.  Thus the $S'$-action is free.  
Let $\mathfrak s':={\rm Lie}\,S'$ (the Lie algebra of $S'$), 
$\widetilde{\mathfrak s}'$ be a maximal solvable subalgebra of $\mathfrak g$ 
containing $\mathfrak s'$ and $\widetilde S'$ be the connected subgroup of $G$ 
with ${\rm Lie}\,\widetilde S'=\widetilde{\mathfrak s}'$.  Since $\mathfrak g$ 
is a real semi-simple Lie algebra and $\widetilde{\mathfrak s}'$ is a maximal 
solvable subalgebra of $\mathfrak g$, $\widetilde{\mathfrak s}'$ contains 
a Cartan subalgebra $\widetilde{\mathfrak a}'$ of $\mathfrak g$.  Let 
$\mathfrak t'$ (resp. $\mathfrak a'$) be the toroidal part (resp. the vector 
part) of $\widetilde{\mathfrak a}'$.  There exists a Cartan decomposition 
$\mathfrak g=\mathfrak f'+\mathfrak p'$ of $\mathfrak g$ with 
$\mathfrak t'\subset\mathfrak f'$ and $\mathfrak a'\subset\mathfrak p'$.  Let 
$\mathfrak g=\mathfrak g'_0+\sum\limits_{\lambda\in\triangle'}
\mathfrak g'_{\lambda}$ be the root space decomposition with respect to 
$\mathfrak a'$ (i.e., $\mathfrak g'_0$ is the centralizer of $\mathfrak a'$ in 
$\mathfrak g$ and $\mathfrak g'_{\lambda}=\{X\in\mathfrak g\,\vert\,{\rm ad}(a)
(X)=\lambda(a)X\,\,{\rm for}\,\,{\rm all}\,\,a\in\mathfrak a'\}$ and 
$\triangle'=\{\lambda\in(\mathfrak a')^{\ast}\setminus \{0\}\,\vert\,
\mathfrak g'_{\lambda}\not=\{0\}\}$).  Let $\mathfrak n':=
\sum\limits_{\lambda\in\triangle'_+}\mathfrak g'_{\lambda}$, where 
$\triangle'_+$ is the positive root system with respect to some lexicographic 
ordering of $\mathfrak a'$.  The algebra $\widetilde{\mathfrak a}'+\mathfrak n'$ is a maximal solvable subalgebra of $\mathfrak g$.  According to a result of 
[Mo], we may assume that $\widetilde{\mathfrak s}'=\widetilde{\mathfrak a}'
+\mathfrak n'$ by retaking $\widetilde{\mathfrak a}'$ if necessary.  
By imitating the proof of Lemma 5.1 of [BT1], it is shown that $\mathfrak a'$ 
is a maximal abelian subspace of $\mathfrak p'$ because 
the $S'$-action has flat section.  
There exists $g\in G$ satisfying ${\rm Ad}(g)(\mathfrak f')=\mathfrak f,\,
{\rm Ad}(g)(\mathfrak p')=\mathfrak p,\,{\rm Ad}(g)(\mathfrak a')
=\mathfrak a$ and ${\rm Ad}(g)({\widetilde{\mathfrak a}}')
=\widetilde{\mathfrak a}$, where ${\rm Ad}$ is the adjoint 
representation of $G$, $\mathfrak a$ and $\widetilde{\mathfrak a}$ are as in 
Introduction.  
Let $\mathfrak s:={\rm Ad}(g)(\mathfrak s')$ and 
$S$ be the connected subgroup of $G$ with ${\rm Lie}\,S=\mathfrak s$.  Since 
the $S$-action is conjugate to the $S'$-action and $S\subset\widetilde AN$, 
we obtain the statement of Theorem A.  \hspace{3.3truecm}q.e.d.

\vspace{0.5truecm}

Let $\mathfrak a$ be a maximal abelian subspace of $\mathfrak p$.  
Fix a lexicographic ordering of $\mathfrak a$.  
Let $\mathfrak g=\mathfrak g_0+\sum\limits_{\lambda\in\triangle}
\mathfrak g_{\lambda}$, $\mathfrak p=\mathfrak a
+\sum\limits_{\lambda\in\triangle_+}\mathfrak p_{\lambda}$ and 
$\mathfrak f=\mathfrak f_0+\sum\limits_{\lambda\in\triangle_+}
\mathfrak f_{\lambda}$ be the root space decompositions of 
$\mathfrak g,\,\mathfrak p$ and $\mathfrak f$ with respect to $\mathfrak a$, 
where we note that 
$$\begin{array}{l}
\displaystyle{\mathfrak g_{\lambda}=\{X\in\mathfrak g\,\vert\,{\rm ad}(a)X
=\lambda(a)X\,\,{\rm for}\,\,{\rm all}\,\,a\in\mathfrak a\}\,\,\,\,
(\lambda\in\triangle),}\\
\displaystyle{\mathfrak p_{\lambda}=\{X\in\mathfrak p\,\vert\,{\rm ad}(a)^2X
=\lambda(a)^2X\,\,{\rm for}\,\,{\rm all}\,\,a\in\mathfrak a\}\,\,\,\,
(\lambda\in\triangle_+),}\\
\displaystyle{\mathfrak f_{\lambda}=\{X\in\mathfrak f\,\vert\,{\rm ad}(a)^2X
=\lambda(a)^2X\,\,{\rm for}\,\,{\rm all}\,\,a\in\mathfrak a\}\,\,\,\,
(\lambda\in\triangle_+\cup\{0\}).}
\end{array}$$
Also, let $\mathfrak g=\mathfrak f+\mathfrak a+\mathfrak n$ be the Iwasawa 
decomposition of $\mathfrak g$ and $G=KAN$ be the corresponding Iwasawa 
decomposition of $G$, where we note that $\mathfrak n
=\sum\limits_{\lambda\in\triangle_+}\mathfrak g_{\lambda}$.  
Now we shall give examples of a solvable group contained in $AN$ whose 
action on $G/K(=AN)$ is complex hyperpolar.  
Denote by $\pi$ the natural projection of $G$ onto $G/K$.  Since $G/K$ is of 
non-compact type, 
$\pi$ gives a diffeomorphism of $AN$ onto $G/K$.  
Denote by $\langle\,\,,\,\,\rangle$ the left-invariant metric of $AN$ induced 
from that of $G/K$ by $\pi\vert_{AN}$.  Also, denote by $\langle\,\,,\,\,
\rangle^G$ the bi-invariant metric of $G$ inducing that of $G/K$.  Note that 
$\langle\,\,,\,\,\rangle\not=\iota^{\ast}\langle\,\,,\,\,\rangle^G$, where 
$\iota$ is the inclusion map of $AN$ into $G$.  
Let ${\it l}$ be a $r$-dimensional subspace of $\mathfrak a+\mathfrak n$ and 
set $\mathfrak s:=(\mathfrak a+\mathfrak n)\ominus{\it l}$, where 
$(\mathfrak a+\mathfrak n)\ominus {\it l}$ denotes the orthogonal complement 
of ${\it l}$ in $\mathfrak a+\mathfrak n$ with respect to 
$\langle\,\,,\,\,\rangle_e$, where $e$ is the identity element of $G$.  
If $\mathfrak s$ is a subalgebra of $\mathfrak a+\mathfrak n$ and 
${\it l}_{\mathfrak p}:={\rm pr}_{\mathfrak p}({\it l})$ 
(${\rm pr}_{\mathfrak p}\,:\,$ the orthogonal projection of $\mathfrak g$ 
onto $\mathfrak p$) is abelian, then the $S$-action ($S:=\exp_G(\mathfrak s)$) 
is a complex hyperpolar action without singular orbit.  
We shall give 
examples of such a subalgebra $\mathfrak s$ of $\mathfrak a+\mathfrak n$ and 
investigate the structure of the $S$-orbit.  

\vspace{0.5truecm}

\noindent
{\it Example 1.} Let $\mathfrak b$ be a $r(\geq1)$-dimensional subspace of 
$\mathfrak a$ and $\mathfrak s_{\mathfrak b}:=(\mathfrak a+\mathfrak n)\ominus
\mathfrak b$.  It is clear that $\mathfrak b_{\mathfrak p}(=\mathfrak b)$ 
is abelian and that $\mathfrak s_{\mathfrak b}$ is 
a subalgebra of $\mathfrak a+\mathfrak n$.  Hence the $S_{\mathfrak b}$-action 
($S_{\mathfrak b}:=\exp_G(\mathfrak s_{\mathfrak b})$) on $G/K$ is a complex 
hyperpolar action without singular orbit.  

\vspace{0.5truecm}

\noindent
{\it Example 2.} Let $\{\lambda_1,\cdots,\lambda_k\}$ be a subset of a simple 
root system $\Pi$ of $\triangle$ such that $H_{\lambda_1},\cdots,
H_{\lambda_k}$ are mutually orthogonal, $\mathfrak b$ be a subspace of 
$\mathfrak a\ominus{\rm Span}\{H_{\lambda_1},\cdots,H_{\lambda_k}\}$ 
(where $\mathfrak b$ may be $\{0\}$) and ${\it l}_i$ ($i=1,\cdots,k$) be 
a one-dimensional subspace of 
${\bf R}H_{\lambda_i}+\mathfrak g_{\lambda_i}$ with 
${\it l}_i\not={\bf R}H_{\lambda_i}$, where $H_{\lambda_i}$ is the 
element of $\mathfrak a$ defined by $\langle H_{\lambda_i},\cdot\rangle=
\lambda_i(\cdot)$ and ${\bf R}H_{\lambda_i}$ is the subspace of $\mathfrak a$ 
spanned by $H_{\lambda_i}$.  Set ${\it l}:=\mathfrak b+\sum\limits_{i=1}^k
{\it l}_i$.  
Then, it follows from Lemma 3.1 (see the below) that ${\it l}_{\mathfrak p}$ 
is abelian and that $\mathfrak s_{\mathfrak b,{\it l}_1,\cdots,{\it l}_k}:=
(\mathfrak a+\mathfrak n)\ominus{\it l}$ is a subalgebra of $\mathfrak a+
\mathfrak n$.  Hence the $S_{\mathfrak b,{\it l}_1,\cdots,{\it l}_k}$-action 
($S_{\mathfrak b,{\it l}_1,\cdots,{\it l}_k}:=\exp_G(
\mathfrak s_{\mathfrak b,{\it l}_1,\cdots,{\it l}_k})$) on $G/K$ is a complex 
hyperpolar action without singular orbit.  

\vspace{0.5truecm}

\noindent
{\bf Lemma 3.1.} {\sl Let ${\it l}$ and 
$\mathfrak s_{\mathfrak b,{\it l}_1,\cdots,{\it l}_k}$ be as in Example 2.  
Then ${\it l}_{\mathfrak p}$ is abelian and 
$\mathfrak s_{\mathfrak b,{\it l}_1,\cdots,{\it l}_k}$ is a subalgebra of 
$\mathfrak a+\mathfrak n$.}

\vspace{0.5truecm}

\noindent
{\it Proof.} 
Let $H\in\mathfrak b$ and $X_i\in{\it l}_i$ ($i=1,\cdots,k$).  Since 
$\lambda_i(H)=0$ and $(X_i)_{\mathfrak p}\in{\bf R}H_{\lambda_i}\oplus
\mathfrak p_{\lambda_i}$, we have $[H,(X_i)_{\mathfrak p}]=0$.  
Fix $i,j\in\{1,\cdots,k\}$ ($i\not=j$).  Since $\lambda_i$ and $\lambda_j$ are 
simple roots and $\langle H_{\lambda_i},H_{\lambda_j}\rangle=0$, 
we have $[(X_i)_{\mathfrak p},(X_j)_{\mathfrak p}]=0$.  
Thus ${\it l}_{\mathfrak p}$ is abelian.  Let 
$V,W\in\mathfrak s_{\mathfrak b,{\it l}_1,\cdots,{\it l}_k}$.  Since 
$\mathfrak s_{\mathfrak b,{\it l}_1,\cdots,{\it l}_k}=(\mathfrak a\ominus
(\mathfrak b+\sum\limits_{i=1}^k{\bf R}H_{\lambda_i}))\oplus
(\sum\limits_{\lambda\in\triangle_+\setminus\{\lambda_1,\cdots,\lambda_k\}}
\mathfrak g_{\lambda})\oplus(\sum\limits_{i=1}^k(({\bf R}H_{\lambda_i}
+\mathfrak g_{\lambda_i})\ominus{\it l}_i))$, $V$ and $W$ are described as 
$V=V_0+\sum\limits_{\lambda\in\triangle_+\setminus\{\lambda_1,\cdots,\lambda_k\}}V_{\lambda}+\sum\limits_{i=1}^kV_i$ and 
$W=W_0+\sum\limits_{\lambda\in\triangle_+\setminus\{\lambda_1,\cdots,\lambda_k\}}W_{\lambda}+\sum\limits_{i=1}^kW_i$, respectively, where $V_0,W_0\in
\mathfrak a\ominus(\mathfrak b+\sum\limits_{i=1}^k{\bf R}H_{\lambda_i}),\,
V_{\lambda},\,W_{\lambda}\in\mathfrak g_{\lambda}$ and $V_i,\,W_i\in
({\bf R}H_{\lambda_i}+\mathfrak g_{\lambda_i})\ominus{\it l}_i$.  Easily 
we have 
$$\begin{array}{l}
\displaystyle{[V,W]\equiv\sum_{\lambda,\mu\in\triangle_+\setminus\{\lambda_1,
\cdots,\lambda_k\}}[V_{\lambda},W_{\mu}]}\\
\hspace{1.5truecm}\displaystyle{
+\sum_{\lambda\in\triangle_+\setminus\{\lambda_1,\cdots,\lambda_k\}}
\sum_{i=1}^k([V_{\lambda},W_i]+[V_i,W_{\lambda}])
+\sum_{i=1}^k\sum_{j=1}^k[V_i,W_j]\quad({\rm mod}\,\,
\mathfrak s_{\mathfrak b,{\it l}_1,\cdots{\it l}_k}).}
\end{array}
$$
Since $\lambda_1,\cdots,\lambda_k$ are simple roots, 
$[V_{\lambda},W_{\mu}],\,[V_{\lambda},W_i],\,[V_i,W_{\lambda}]$ 
and $[V_i,W_j]$ ($\lambda,\mu\in
\triangle_+\setminus\{\lambda_1,\cdots,\lambda_k\},\,1\leq i,j\leq k$) 
belong to $S_{\mathfrak b,{\it l}_1,\cdots,{\it l}_k}$.  
Therefore we have $[V,W]\in\mathfrak s_{\mathfrak b,{\it l}_1,\cdots,
{\mathfrak l}_k}$.  Thus $\mathfrak s_{\mathfrak b,{\it l}_1,\cdots{\it l}_k}$ 
is a subalgebra of $\mathfrak a+\mathfrak n$.  \hspace{8.4truecm} q.e.d.

\vspace{0.5truecm}

For the orbit $S_{\mathfrak b,{\it l}_1,\cdots,{\it l}_k}\cdot e$, we have 
the following facts.  

\vspace{0.5truecm}

\noindent
{\bf Lemma 3.2.} {\sl Let 
$\mathfrak s_{\mathfrak b,{\it l}_1,\cdots,{\it l}_k}$ be as in Example 2, 
$\xi_0\in\mathfrak b,\,\xi^i_{t_i}:=
\frac{1}{\cosh(\vert\lambda_i\vert t_i)}\xi^i
-\frac{1}{\vert\lambda_i\vert}\tanh(\vert\lambda_i\vert t_i)H_{\lambda_i}$ be 
a unit vector of ${\it l}_i$ ($i=1,\cdots,k$), where $\xi^i$ is a unit vector 
of $\mathfrak g_{\lambda_i}$.  Denote by $A$ the shape tensor of the orbit 
$S_{\mathfrak b,{\it l}_1,\cdots,{\it l}_k}\cdot e\,(\subset AN)$.  
Then, for $A_{\xi_0}$ and $A_{\xi^i_{t_i}}$, 
the following statements ${\rm(i)}\sim{\rm(vii)}$ hold:

{\rm(i)} For $X\in\mathfrak a\ominus(\mathfrak b+\sum\limits_{i=1}^k{\rm R}
H_{\lambda_i})$, we have $A_{\xi_0}X=A_{\xi^i_{t_i}}X=0$ 
($i=1,\cdots,k$).

{\rm (ii)} For $X\in{\rm Ker}({\rm ad}(\xi^i)\vert_{\mathfrak g_{\lambda_i}})
\ominus{\bf R}\xi^i$, 
we have $A_{\xi_0}X=0$ and $A_{\xi^i_{t_i}}X=-\vert\lambda_i\vert
\tanh(\vert\lambda_i\vert t_i)X$.  

{\rm (iii)} Assume that $2\lambda_i\in\triangle_+$.  For 
$X\in\mathfrak g_{2\lambda_i}$, we have 
$A_{\xi_0}([\theta\xi^i,X])=0$ and 
$$\begin{array}{l}
\displaystyle{A_{\xi^i_{t_i}}X=-2\vert\lambda_i\vert\tanh(
\vert\lambda_i\vert t_i)X-\frac{1}{2\cosh(\vert\lambda_i\vert t_i)}
[\theta\xi^i,X],}\\
\displaystyle{A_{\xi^i_{t_i}}([\theta\xi^i,X])
=-\frac{\vert\lambda_i\vert^2}{\cosh(\vert\lambda_i\vert t_i)}X
-\vert\lambda_i\vert\tanh(\vert\lambda_i\vert t_i)[\theta\xi^i,X],}
\end{array}$$
where $\theta$ is the Cartan involution of $\mathfrak g$ with 
${\rm Fix}\,\theta=\mathfrak f$.  

{\rm(iv)} For $X\in({\bf R}\xi^i+{\bf R}H_{\lambda_i})\ominus{\it l}_i$, 
we have $A_{\xi_0}X=0$ and $A_{\xi^i_{t_i}}X=-\vert\lambda_i\vert
\tanh(\vert\lambda_i\vert t_i)X$.

{\rm (v)} For $X\in(\mathfrak g_{\lambda_j}\ominus{\bf R}\xi^j)
+(({\bf R}\xi^j+{\bf R}H_{\lambda_j})\ominus{\it l}_j)
+\mathfrak g_{2\lambda_j}$ ($j\not=i$), we have 
$A_{\xi_0}X=A_{\xi^i_{t_i}}X=0$.

{\rm(vi)} For $X\in\mathfrak g_{\mu}$ ($\mu\in\triangle_+\setminus
\{\lambda_1,\cdots,\lambda_k\}$), we have $A_{\xi_0}X=\mu(\xi_0)X$.

{\rm (vii)} Let 
$\displaystyle{k_i:=\exp\left(\frac{\pi}{\sqrt2\vert\lambda_i\vert}
(\xi^i+\theta\xi^i)\right)}$, where $\exp$ is the exponential map of $G$.  
Then ${\rm Ad}(k_i)\circ A_{\xi^i_{t_i}}=-A_{\xi^i_{t_i}}\circ{\rm Ad}(k_i)$ 
holds over $\mathfrak n\ominus\sum\limits_{i=1}^k(\mathfrak g_{\lambda_i}
+\mathfrak g_{2\lambda_i})$, where ${\rm Ad}$ is the adjoint representation 
of $G$.}

\vspace{0.5truecm}

\noindent
{\it Proof.} Let ${\rm pr}^1_{\mathfrak a+\mathfrak n}$ (resp. 
${\rm pr}^2_{\mathfrak a+\mathfrak n}$) be the projection of $\mathfrak g$ 
onto $\mathfrak a+\mathfrak n$ with respect to the decomposition 
$\mathfrak g=\mathfrak f+(\mathfrak a+\mathfrak n)$ (resp. 
$\mathfrak g=(\mathfrak f_0+\sum\limits_{\lambda\in\triangle_+}
\mathfrak p_{\lambda})+(\mathfrak a+\mathfrak n)$), ${\rm pr}_{\mathfrak f}$ 
(resp. ${\rm pr}_{\mathfrak p}$) be the projection of $\mathfrak g$ onto 
$\mathfrak f$ (resp. $\mathfrak p$) with respect to the decomposition 
$\mathfrak g=\mathfrak f+\mathfrak p$ and ${\rm pr}_{\mathfrak f_0}$ be 
the projection of $\mathfrak g$ onto $\mathfrak f_0$ with respect to the 
decomposition $\mathfrak g=\mathfrak f_0+(\mathfrak a+
\sum\limits_{\lambda\in\triangle}\mathfrak g_{\lambda})$.  Then we have 
$${\rm pr}_{\mathfrak p}\circ{\rm pr}^1_{\mathfrak a+\mathfrak n}
={\rm pr}_{\mathfrak p}\,\,\,\,{\rm and}\,\,\,\,{\rm pr}_{\mathfrak f}\circ
{\rm pr}^2_{\mathfrak a+\mathfrak n}={\rm pr}_{\mathfrak f}-
{\rm pr}_{\mathfrak f_0}.\leqno{(3.1)}$$
Let $H\in\mathfrak a,\,N_1,N_2\in\mathfrak n$ 
and $E\in\mathfrak g_{\lambda}$ ($\lambda\in\triangle_+$).  Denote by 
${\rm ad}(H)^{\ast}$ (resp. ${\rm ad}(E)^{\ast}$) the adjoint operator of 
${\rm ad}(H)$ (resp. ${\rm ad}(E)$) $:\mathfrak a+\mathfrak n\to\mathfrak a
+\mathfrak n$ with respect to $\langle\,\,,\,\,\rangle_e$.  
Easily we can show 
$${\rm ad}(H)^{\ast}={\rm ad}(H).\leqno{(3.2)}$$
For simplicity, we denote ${\rm pr}_{\mathfrak f}(\cdot)$ (resp. 
${\rm pr}_{\mathfrak p}(\cdot)$) by $(\cdot)_{\mathfrak f}$ (resp. 
$(\cdot)_{\mathfrak p}$).  From $(3.1)$ and the skew-symmetricness of 
${\rm ad}(\cdot)$ with respect to $\langle\,\,,\,\,\rangle^G_e$, we have 
$$\begin{array}{l}
\displaystyle{\langle{\rm ad}(E)N_1,N_2\rangle_e=
\langle{\rm ad}(E_{\mathfrak f})((N_1)_{\mathfrak p})+{\rm ad}(E_{\mathfrak p})
(((N_1)_{\mathfrak f}),\,(N_2)_{\mathfrak p}\rangle^G_e}\\
\hspace{2.7truecm}\displaystyle{=
-\langle(N_1)_{\mathfrak p},\,{\rm ad}(E_{\mathfrak f})((N_2)_{\mathfrak p})
\rangle^G_e-\langle(N_1)_{\mathfrak f},\,{\rm ad}(E_{\mathfrak p})
((N_2)_{\mathfrak p})\rangle^G_e}\\
\hspace{2.7truecm}\displaystyle{=
-\langle(N_1)_{\mathfrak p},\,({\rm pr}^1_{\mathfrak a+\mathfrak n}
({\rm ad}(E_{\mathfrak f})N_2))_{\mathfrak p}\rangle^G_e}\\
\hspace{3.2truecm}\displaystyle{
-\langle(N_1)_{\mathfrak f},\,({\rm pr}^2_{\mathfrak a+\mathfrak n}
({\rm ad}(E_{\mathfrak p})N_2))_{\mathfrak f}+{\rm pr}_{\mathfrak f_0}
({\rm ad}(E_{\mathfrak p})N_2)\rangle^G_e}\\
\hspace{2.7truecm}\displaystyle{=
-\langle N_1,\,{\rm pr}^1_{\mathfrak a+\mathfrak n}
({\rm ad}(E_{\mathfrak f})N_2)\rangle_e
+\langle N_1,\,{\rm pr}^2_{\mathfrak a+\mathfrak n}
({\rm ad}(E_{\mathfrak p})N_2)\rangle_e}
\end{array}$$
and hence 
$${\rm pr}_{\mathfrak n}({\rm ad}(E)^{\ast}N_2)={\rm pr}_{\mathfrak n}
(-{\rm pr}^1_{\mathfrak a+\mathfrak n}({\rm ad}(E_{\mathfrak f})N_2)
+{\rm pr}^2_{\mathfrak a+\mathfrak n}({\rm ad}(E_{\mathfrak p})N_2)),$$
where ${\rm pr}_{\mathfrak n}$ is the projection of $\mathfrak a+\mathfrak n$ 
onto $\mathfrak n$.  Also, we have 
$$\langle{\rm ad}(E)H,N_2\rangle_e=-\lambda(H)\langle E,N_2\rangle_e
=-\langle H,\langle E,N_2\rangle_eH_{\lambda}\rangle_e$$
and hence ${\rm pr}_{\mathfrak a}({\rm ad}(E)^{\ast}N_2)=-\langle E,N_2
\rangle_eH_{\lambda}$, where ${\rm pr}_{\mathfrak a}$ is the projection 
of $\mathfrak a+\mathfrak n$ onto $\mathfrak a$.  Also, we can show 
${\rm ad}(E)^{\ast}H=0$.  Therefore, we have 
$${\rm ad}(E)^{\ast}=\left\{
\begin{array}{cl}
\displaystyle{0} & \displaystyle{{\rm on}\,\,\mathfrak a}\\
\displaystyle{
\begin{array}{l}
\displaystyle{-\langle E,\cdot\rangle_e\otimes H_{\lambda}-
{\rm pr}_{\mathfrak n}\circ{\rm pr}^1_{\mathfrak a+\mathfrak n}\circ
{\rm ad}(E_{\mathfrak f})}\\
\displaystyle{+{\rm pr}_{\mathfrak n}\circ{\rm pr}^2_{\mathfrak a+\mathfrak n}
\circ{\rm ad}(E_{\mathfrak p})}
\end{array}
} & \displaystyle{{\rm on}\,\,\mathfrak n}
\end{array}
\right.
\leqno{(3.3)}$$
On the other hand, according to the Koszul's formula, we have 
$$\begin{array}{l}
\displaystyle{\langle A_{\xi}X,Y\rangle_e=\frac12\left(
\langle[X,Y],\xi\rangle_e-\langle[Y,\xi],X\rangle_e
+\langle[\xi,X],Y\rangle_e\right)}\\
\hspace{2.1truecm}\displaystyle{=\frac12\langle({\rm ad}(\xi)+{\rm ad}(\xi)^{\ast})X,Y\rangle_e}
\end{array}$$
for any $X,Y\in T_e(S_{\mathfrak b,{\it l}_1,\cdots,{\it l}_k}\cdot e)=
\mathfrak s_{\mathfrak b,{\it l}_1,\cdots,{\it l}_k}$ and any $\xi\in 
T^{\perp}_e(S_{\mathfrak b,{\it l}_1,\cdots,{\it l}_k}\cdot e)
=\mathfrak b+\sum\limits_{i=1}^k{\it l}_i$.  That is, we have 
$$A_{\xi}=\frac12{\rm pr}_T\circ({\rm ad}(\xi)+{\rm ad}(\xi)^{\ast}),
\leqno{(3.4)}$$
where ${\rm pr}_T$ is the orthogonal projection of $\mathfrak a+\mathfrak n$ 
onto $\mathfrak s_{\mathfrak b,{\it l}_1,\cdots,{\it l}_k}$.  From $(3.2)$ and 
$(3.4)$, we have 
$$A_{\xi_0}X=\left\{
\begin{array}{cl}
\displaystyle{0}&\displaystyle{(X\in\mathfrak s_{\mathfrak b,{\it l}_1,
\cdots,{\it l}_k}\ominus\sum_{\mu\in\triangle_+\setminus\{\lambda_1,\cdots,
\lambda_k\}}\mathfrak g_{\lambda})}\\
\displaystyle{\mu(\xi_0)X}&\displaystyle{(X\in\mathfrak g_{\mu}),}
\end{array}\right.$$
where $\mu\in\triangle_+\setminus\{\lambda_1,\cdots,\lambda_k\}$.  From 
$(3.3)$ and $(3.4)$, we have 
$$A_{\xi^i_{t_i}}X=0\qquad(X\in\mathfrak a\ominus
(\mathfrak b+\sum_{i=1}^k{\bf R}H_{\lambda_i})).$$
Set $\mathfrak g^K_{\lambda_j}:={\rm Ker}({\rm ad}(\xi^j)
\vert_{\mathfrak g_{\lambda_j}})$ and 
$\mathfrak g_{\lambda_j}^I:={\rm Im}({\rm ad}(\theta\xi^j)
\vert_{\mathfrak g_{2\lambda_j}})$ ($j=1,\cdots,k$).  Then we have 
$\mathfrak g_{\lambda_j}=\mathfrak g^K_{\lambda_j}\oplus
\mathfrak g^I_{\lambda_j}$.  By simple calculations, it is shown that this 
decomposition is orthogonal with respect to $\langle\,\,,\,\,\rangle_e$.  
If $X\in\mathfrak g^K_{\lambda_j}\ominus{\bf R}\xi^j$, then it follows from 
$(3.2),\,(3.3),\,(3.4),\,\lambda_i,\lambda_j\in\Pi$ and 
$\langle H_{\lambda_i},H_{\lambda_j}\rangle=0$ (when $i\not=j$) that 
$$A_{\xi^i_{t_i}}X=\left\{
\begin{array}{cc}
\displaystyle{-\vert\lambda_i\vert\tanh(\vert\lambda_i\vert t_i)X} & 
\displaystyle{(i=j)}\\
\displaystyle{0} & \displaystyle{(i\not=j).}
\end{array}
\right.$$
If $X\in\mathfrak g_{2\lambda_j}$, then it follows from $(3.2),\,(3.3),\,
(3.4),\,\lambda_i,\lambda_j\in\Pi$ and 
$\langle H_{\lambda_i},H_{\lambda_j}\rangle=0$ (when $i\not=j$) that 
$$A_{\xi^i_{t_i}}X=
\left\{
\begin{array}{cc}
\displaystyle{-2\vert\lambda_i\vert\tanh(\vert\lambda_i\vert t_i)X
-\frac{1}{2\cosh(\vert\lambda_i\vert t_i)}[\theta\xi^i,X]} & 
\displaystyle{(i=j)}\\
\displaystyle{0} & \displaystyle{(i\not=j).}
\end{array}
\right.$$
Also, we have 
$$A_{\xi^i_{t_i}}([\theta\xi^i,X])=
-\frac{\vert\lambda_i\vert^2}
{\cosh(\vert\lambda_i\vert t_i)}X-\vert\lambda_i\vert
\tanh(\vert\lambda_i\vert t_i)[\theta\xi^i,X].$$
Let $X:=\tanh(\vert\lambda_j\vert t_j)\xi^j+\frac{1}
{\vert\lambda_j\vert\cosh(\vert\lambda_j\vert t_j)}H_{\lambda_j}$, which is 
a unit vector of $({\bf R}\xi^j+{\bf R}H_{\lambda_j})\ominus{\it l}_j$.  
From $(3.2),\,(3.3),\,(3.4),\,\lambda_i,\lambda_j\in\Pi$ and 
$\langle H_{\lambda_i},H_{\lambda_j}\rangle=0$ (when $i\not=j$), we have 
$$\begin{array}{l}
\displaystyle{A_{\xi^i_{t_i}}X=-\frac12\vert\lambda_i\vert
\tanh(\vert\lambda_i\vert t_i)X+\frac{1}{2\cosh(\vert\lambda_i\vert t_i)}
{\rm pr}_T({\rm ad}(\xi^i)^{\ast}X)}\\
\hspace{2truecm}\displaystyle{-\frac{1}{2\vert\lambda_i\vert}
\tanh(\vert\lambda_i\vert t_i){\rm pr}_T({\rm ad}(H_{\lambda})^{\ast}X)}\\
\hspace{1.2truecm}\displaystyle{=
\left\{
\begin{array}{cc}
\displaystyle{-\vert\lambda_i\vert\tanh(\vert\lambda_i\vert t_i)X} & 
\displaystyle{(i=j)}\\
\displaystyle{0} & \displaystyle{(i\not=j).}
\end{array}\right.}
\end{array}$$
This completes the proof of ${\rm(i)}\sim{\rm(vi)}$.  
Finally we shall show the statement (vii).  
Let $X\in\mathfrak n\ominus\sum\limits_{i=1}^k(\mathfrak g_{\lambda_i}
+\mathfrak g_{2\lambda_i})$ and $k_i$ be as in the statement (vii).  From 
$(3.2),\,(3.3),\,(3.4),\,\lambda_j\in\Pi$ ($j=1,\cdots,k$) and 
$\langle H_{\lambda_i},H_{\lambda_j}\rangle=0$ (when $i\not=j$), we have 
$$A_{\xi^i_{t_i}}X=\frac{1}{\cosh(\vert\lambda_i\vert t_i)}
[\xi^i_{\mathfrak p},X]-\frac{1}{\vert\lambda_i\vert}
\tanh(\vert\lambda_i\vert t_i)[H_{\lambda_i},X].$$
By operating ${\rm Ad}(k_i)$ to both sides of this relation, we have 
$${\rm Ad}(k_i)(A_{\xi^i_{t_i}}X)=-A_{\xi^i_{t_i}}({\rm Ad}(k_i)X),$$
where we use ${\rm Ad}(k_i)(\xi^i_{\mathfrak p})=-\xi^i_{\mathfrak p}$ and 
${\rm Ad}(k_i)(H_{\lambda_i})=-H_{\lambda_i}$.  Thus the statement (vii) 
is shown.  \hspace{12.7truecm}q.e.d.

\vspace{0.5truecm}

Also, we have the following fact.  

\vspace{0.5truecm}

\noindent
{\bf Lemma 3.3.} {\sl Let $\mathfrak s_{\mathfrak b,{\it l}_1,\cdots,
{\it l}_k}$ be as in Example 2 and $\bar{\it l}_i$ be the orthogonal 
projection of ${\it l}_i$ onto $\mathfrak g_{\lambda_i}$.  Set 
$\mathfrak s_{\mathfrak b,\bar{\it l}_1,\cdots,\bar{\it l}_k}:=
(\mathfrak a+\mathfrak n)\ominus(\mathfrak b+\sum\limits_{i=1}^k
\bar{\it l}_i)$ and $S_{\mathfrak b,\bar{\it l}_1,\cdots,\bar{\it l}_k}:=
\exp_G(\mathfrak s_{\mathfrak b,\bar{\it l}_1,\cdots,\bar{\it l}_k})$.  Then 
the $S_{\mathfrak b,\bar{\it l}_1,\cdots,\bar{\it l}_k}$-action is conjugate 
to the $S_{\mathfrak b,{\it l}_1,\cdots,{\it l}_k}$-action.}

\vspace{0.5truecm}

\noindent
{\it Proof.} Denote by $\nabla$ the Levi-Civita connection of the 
left-invariant metric of $AN$.  
Let $H$ be a vector of $\mathfrak b$, $\xi^i$ be a unit vector of 
$\bar{\it l}_i$ ($i=1,\cdots,k$) and $\gamma_{\xi^i}$ be the geodesic in $AN$ 
with $\dot{\gamma}_{\xi^i}(0)=\xi^i$.  
Let $t_i$ be a real number with 
$\frac{1}{\cosh(\vert\lambda_i\vert t_i)}\xi^i-\tanh(\vert\lambda_i\vert t_i)
H_{\lambda_i}\in{\it l}_i$ ($i=1,\cdots,k$).  
Denote by the same symbols $H,\,\xi^i$ and $H_{\lambda_i}$ the left-invariant 
vector fields arising from $H,\,\xi^i$ and $H_{\lambda_i}$, respectively.  
By using the relation $(5.4)$ of Section 5 of [Mi] (arising the Koszul formula 
for the left-invariant vector fields), we can show 
$$\begin{array}{l}
\displaystyle{\nabla_{\xi^1}\xi^1=\vert\lambda_1\vert H_{\lambda_1},\,
\nabla_{\xi^1}H_{\lambda_1}=-\vert\lambda_1\vert\xi^1}\\
\displaystyle{\nabla_{\xi^1}\xi^i=\nabla_{\xi^1}H=
\nabla_{H_{\lambda_1}}\xi^1=\nabla_{H_{\lambda_1}}\xi^i
=\nabla_{H_{\lambda_1}}H_{\lambda_1}=\nabla_{H_{\lambda_1}}H=0,}
\end{array}$$
where $i=2,\cdots,k$.  From 
$\nabla_{\xi^1}\xi^1=\vert\lambda_1\vert H_{\lambda_1},\,
\nabla_{\xi^1}H_{\lambda_1}=-\vert\lambda_1\vert\xi^1,\,
\nabla_{H_{\lambda_1}}\xi^1=\nabla_{H_{\lambda_1}}H_{\lambda_1}=0$, 
it follows that $\exp{\bf R}\{\xi^1,H_{\lambda_1}\}$ is a totally geodesic 
subgroup of $AN$.  Hence $\dot{\gamma}_{\xi^1}(t)$ is expressed as 
$\dot{\gamma}_{\xi^1}(t)=a(t)(H_{\lambda_1})_{\gamma_{\xi^1}(t)}
+b(t)(\xi^1)_{\gamma_{\xi^1}(t)}$.  Furthermore, we have 
$\nabla_{\dot{\gamma}_{\xi^1}}\dot{\gamma}_{\xi^1}
=(a'+\vert\lambda_1\vert b^2)H_{\lambda_1}+(b'-\vert\lambda_1\vert ab)
\xi^1=0$, that is, 
$a'=-\vert\lambda_1\vert b^2$ and $b'=\vert\lambda_1\vert ab$.  
By solving this differential equation under the initial conditions $a(0)=0$ 
and $b(0)=1$, we have $a(t)=-\tanh(\vert\lambda_1\vert t)$ 
and $b(t)=\frac{1}{\cosh(\vert\lambda_1\vert t)}$.  
Hence we obtain 
$\dot{\gamma}_{\xi^1}(t)=
\frac{1}{\cosh(\vert\lambda_1\vert t)}(\xi^1)_{\gamma_{\xi^1}(t)}
-\tanh(\vert\lambda_1\vert t)(H_{\lambda_1})_{\gamma_{\xi^1}(t)}$.  
From $\nabla_{\xi^1}\xi^i=\nabla_{\xi^1}H=\nabla_{H_{\lambda_1}}\xi^i=
\nabla_{H_{\lambda_1}}H=0$ ($i=2,\cdots,k$), it follows that 
$\xi^i$ ($i=2,\cdots,k$) and $H$ are parallel along $\gamma_{\xi^1}$ (with 
respect to $\nabla$).  Denote by $P_{\gamma_{\xi^1}\vert_{[0,t]}}$ the 
parallel translation along $\gamma_{\xi^1}\vert_{[0,t]}$ (with 
respect to $\nabla$) and $L_{\gamma_{\xi^1}(t)}$ the left translation by 
$\gamma_{\xi^1}(t)$.  From the above facts, we have 
$$
\begin{array}{l}
\displaystyle{
T^{\perp}_{\gamma_{\xi^1}(t_1)}(S_{\mathfrak b,\bar{\it l}_1,\cdots,
\bar{\it l}_k})=P_{\gamma_{\xi^1}\vert_{[0,t_1]}}(\mathfrak b
+\sum_{i=1}^k\bar{\it l}_i)
=(L_{\gamma_{\xi^1}(t_1)})_{\ast}(\mathfrak b
+\sum_{i=2}^k\bar{\it l}_i+{\it l}_1)}\\
\hspace{2.7truecm}\displaystyle{=(L_{\gamma_{\xi^1}(t_1)})_{\ast}
(T^{\perp}_eS_{\mathfrak b,{\it l}_1,\bar{\it l}_2,\cdots,\bar{\it l}_k}),}
\end{array}$$
which implies $\gamma_{\xi^1}(t_1)^{-1}S_{\mathfrak b,\bar{\it l}_1,
\cdots,\bar{\it l}_k}\gamma_{\xi^1}(t_1)=S_{\mathfrak b,{\it l}_1,
\bar{\it l}_2,\cdots,\bar{\it l}_k}$.  By repeating the same discussion, we 
obtain 
$$(\gamma_{\xi^1}(t_1)\cdots\gamma_{\xi^k}(t_k))^{-1}
S_{\mathfrak b,\bar{\it l}_1,\cdots,\bar{\it l}_k}(\gamma_{\xi^1}(t_1)\cdots
\gamma_{\xi^k}(t_k))=S_{\mathfrak b,{\it l}_1,\cdots,{\it l}_k}.$$
Thus the $S_{\mathfrak b,\bar{\it l}_1,\cdots,\bar{\it l}_k}$-action is 
conjugate to the $S_{\mathfrak b,{\it l}_1,\cdots,{\it l}_k}$-action.  
\hspace{4.1truecm}q.e.d.

\vspace{0.5truecm}

For parallel submanifolds of a proper complex equifocal submanifold and 
a curvature-adapted complex equifocal submanifold, we have the following 
facts.  

\vspace{0.5truecm}

\noindent
{\bf Lemma 3.4.} {\sl {\rm (i)} All parallel submanifolds of a proper 
complex equifocal submanifold are proper complex equifocal.  

{\rm (ii)} All parallel submanifolds of a curvature-adapted complex equifocal 
submanifold are curvature-adapted and complex equifocal.}

\vspace{0.5truecm}

\noindent
{\it Proof.} First we shall show the statement (i).  Let $M$ be a proper 
complex equifocal submanifold in a symmetric space $G/K$ of non-compact type 
and $\widetilde v$ be the parallel normal vector field of $M$ which is not 
a focal normal vector field.  Denote by $\eta_{\widetilde v}$ the end-point 
map for $\widetilde v$ and $M_{\widetilde v}:=\eta_{\widetilde v}(M)$, which 
is a parallel submanifold of $M$.  The vector field $\widetilde v$ is regarded 
as a parallel normal vector field of the complexification $M^{\bf c}$ 
along $M$.  
Let ${\widetilde v}^L$ be the horizontal lift of $\widetilde v$ to 
$H^0([0,1],\mathfrak g^{\bf c})$ by the anti-Kaehlerian submersion 
$\pi^{\bf c}\circ\phi^{\bf c}:H^0([0,1],\mathfrak g^{\bf c})\to 
G^{\bf c}/K^{\bf c}$, which is a parallel normal vector field of 
${\widetilde M}^{\bf c}(:=(\pi^{\bf c}\circ\phi^{\bf c})^{-1}(M^{\bf c}))$.  
Set $\widetilde{M^{\bf c}}_{{\widetilde v}^L}:=\eta_{{\widetilde v}^L}
(\widetilde{M^{\bf c}})$, where $\eta_{{\widetilde v}^L}$ is the end-point map 
for ${\widetilde v}^L$.  Note that $\widetilde{M^{\bf c}}_{{\widetilde v}^L}
=(\pi^{\bf c}\circ\phi^{\bf c})^{-1}((M_{\widetilde v})^{\bf c})$.  Denote by 
$\widetilde A$ and $\widetilde A^{{\widetilde v}^L}$ the shape tensors of 
$\widetilde{M^{\bf c}}$ and $\widetilde{M^{\bf c}}_{{\widetilde v}^L}$, 
respectively.  Let $\{\lambda_i\,\vert\,i\in I\}$ be the set of all complex 
principal curvatures of $\widetilde{M^{\bf c}}$ and $E_i$ be the complex 
curvature distribution for $\lambda_i$.  Then, according to Lemma 3.2 of 
[Koi4], we have 
$$\widetilde A_w^{{\widetilde v}^L}\vert_{(E_i)_u}
=\frac{(\lambda_i)_u(w)}{1-(\lambda_i)_u({\widetilde v}^L_u)}{\rm id}\,\,\,\,
(i\in I,\,\,u\in\widetilde{M^{\bf c}}_{{\widetilde v}^L}),\leqno{(3.5)}$$
where we note that $T_{\eta_{{\widetilde v}^L}(u)}
\widetilde{M^{\bf c}}_{{\widetilde v}^L}=T_u\widetilde{M^{\bf c}}(=
\displaystyle{\overline{\mathop{\oplus}_{i\in I}(E_i)_u})}$.  This implies 
that $\widetilde{M^{\bf c}}_{{\widetilde v}^L}$ is proper anti-Kaehlerian 
isoparametric, that is, $M_{\widetilde v}$ is proper complex equifocal.  Thus 
the statement (i) is shown.  Next we shall show the statement (ii).  Let $M$ 
be a curvature-adapted complex equifocal submanifold in $G/K$ and 
$\widetilde v$ be the parallel normal vector field of $M$.  Set 
$M_{\widetilde v}:=\eta_{\widetilde v}(M)$.  Denote by $A$ and 
$A^{\widetilde v}$ the shape tensors of $M$ and $M_{\widetilde v}$, 
repsectively.  
Let $w\in T^{\perp}_xM$.  Without loss of generality, we may assume that 
$x=eK$.  Let $\mathfrak a$ be a maximal abelian subspace of 
$\mathfrak p:=T_{eK}(G/K)$ containing $T^{\perp}_{eK}M$ and 
$\mathfrak p=\mathfrak a+\sum\limits_{\alpha\in\triangle_+}
\mathfrak p_{\alpha}$ be the root space decomposition with respect to 
$\mathfrak a$.  Let $X\in{\rm Ker}(A_v-\lambda\,{\rm id})\cap{\rm Ker}
(A_w-\mu\,{\rm id})\cap\mathfrak p_{\alpha}$ ($\lambda\in{\rm Spec}\,A_v,\,\,
\mu\in{\rm Spec}\,A_w,\,\,\alpha\in\triangle_+$).  Let $\widetilde w$ be the 
parallel tangent vector field on the (flat) section $\Sigma$ of $M$ through 
$eK$ with $\widetilde w_{eK}=w$.  Since $M_{\widetilde v}$ is regarded as a 
partial tube over $M$, it follows from (ii) of Corollary 3.2 in [Koi3] that 
$$\begin{array}{l}
\displaystyle{
(A^{\widetilde v})_{\widetilde w_{\eta_{\widetilde v}(eK)}}
((\eta_{\widetilde v})_{\ast}X)=
\frac{1}{\alpha(v)-\lambda\tanh\,\alpha(v)}
\{-\alpha(v)\alpha(w)\tanh\,\alpha(v)}\\
\hspace{3.1truecm}\displaystyle{+\lambda\left(
1-\frac{\tanh\,\alpha(v)}{\alpha(v)}\right)\alpha(w)
+\mu\tanh\,\alpha(v)\}(\eta_{\widetilde v})_{\ast}X.}
\end{array}\leqno{(3.6)}$$
Let $Z$ be the element of $\mathfrak p$ with $\exp_G(Z)K=\eta_{\widetilde v}
(eK)$.  For simplicity, set $g:=\exp_G(Z)$.  Since $g_{\ast}:\mathfrak p\to 
T_{\eta_{\widetilde v}(eK)}(G/K)$ is the parallel trnaslation along the normal 
geodesic $\gamma_Z(\displaystyle{\mathop{\Leftrightarrow}_{{\rm def}}\,\,}$
\newline
$\gamma_Z(t):=\exp_G(tZ)K)$, it follows from $(3.1)$ of [Koi3] that 
$$\begin{array}{l}
\displaystyle{(\eta_{\widetilde v})_{\ast}X=g_{\ast}
(D^{co}_v(X)-D^{si}_v(A_vX))}\\
\hspace{1.4truecm}\displaystyle{=\left(\cosh\,\alpha(v)-\lambda
\frac{\sinh\,\alpha(v)}{\alpha(v)}\right)g_{\ast}X\in 
g_{\ast}\mathfrak p_{\alpha}.}
\end{array}$$
Also, we have $g_{\ast}^{-1}(T^{\perp}_{\eta_{\widetilde v}(eK)}
M_{\widetilde v})=T^{\perp}_{eK}M\subset\mathfrak a$.  Hence we have 
$R((\eta_{\widetilde v})_{\ast}X,\widetilde w_{\eta_{\widetilde v}(eK)})
\widetilde w_{\eta_{\widetilde v}(eK)}=-\alpha(w)^2
(\eta_{\widetilde v})_{\ast}X$, which together with $(3.6)$ implies 
$$[(A^{\widetilde v})_{\widetilde w_{\eta_{\widetilde v}(eK)}},
R(\cdot,\widetilde w_{\eta_{\widetilde v}(eK)})
\widetilde w_{\eta_{\widetilde v}(eK)}]((\eta_{\widetilde v})_{\ast}X)=0.$$
Therefore, it follows from the arbitrariness of $X$ that 
$[(A^{\widetilde v})_{\widetilde w_{\eta_{\widetilde v}(eK)}},
R(\cdot,\widetilde w_{\eta_{\widetilde v}(eK)})
\widetilde w_{\eta_{\widetilde v}(eK)}]$ vanishes over 
$(\eta_{\widetilde v})_{\ast}({\rm Ker}(A_v-\lambda\,{\rm id})\cap
{\rm Ker}(A_w-\mu\,{\rm id})\cap\mathfrak p_{\alpha})$.  Since $M$ is 
curvature-adapted, we have 
$$\mathop{\oplus}_{\lambda\in{\rm Spec}\,A_v}
\mathop{\oplus}_{\mu\in{\rm Spec}\,A_w}
\mathop{\oplus}_{\alpha\in\triangle_+}(\eta_{\widetilde v})_{\ast}
({\rm Ker}(A_v-\lambda\,{\rm id})\cap {\rm Ker}(A_w-\mu\,{\rm id})\cap
\mathfrak p_{\alpha})=T_{\eta_{\widetilde v}(eK)}M_{\widetilde v}.$$
Hence we have $[(A^{\widetilde v})_{\widetilde w_{\eta_{\widetilde v}(eK)}},
R(\cdot,\widetilde w_{\eta_{\widetilde v}(eK)})
\widetilde w_{\eta_{\widetilde v}(eK)}]=0$.  
Therefore, it follows from the arbitrariness of $w$ that $M_{\widetilde v}$ is 
curvature-adapted.  It is clear that $M_{\widetilde v}$ is complex equifocal.  
Thus the statement (ii) is shown.  \hspace{8.4truecm}q.e.d.

\vspace{0.5truecm}

For the $S_{\mathfrak b}$-action and the $S_{\mathfrak b,{\it l}_1,\cdots,
{\it l}_k}$-action, we have the following facts.  

\vspace{0.5truecm}

\noindent
{\bf Proposition 3.5.} {\sl {\rm (i)} All orbits of the $S_{\mathfrak b}$-action are curvature-adapted but they are not proper complex equifocal.  

{\rm (ii)} Let $\lambda_1,\cdots,\lambda_k\,(\in\triangle_+)$ be as in 
Example 2.  If the root system $\triangle$ of $G/K$ is non-reduced and 
$2\lambda_{i_0}\in\triangle_+$ for some $i_0\in\{1,\cdots,k\}$, then all 
orbits of the $S_{\mathfrak b,{\it l}_1,\cdots,{\it l}_k}$-action are not 
curvature-adapted.  
Also, if $\mathfrak b\not=\{0\}$, then they are not proper complex equifocal.}

\vspace{0.5truecm}

\noindent
{\it Proof.} First we shall show the statement (i).  
The group $S_{\mathfrak b}$ acts isometrically on 
$(AN,\langle\,\,,\,\,\rangle)$.  
Denote by $A$ the shape tensor of the orbit $S_{\mathfrak b}\cdot e$ in $AN$.  
Since 
$\langle\,\,,\,\,\rangle$ is left-invariant, it follows from the Koszul 
formula that $\langle A_vX,Y\rangle=\langle{\rm ad}(v)X,Y\rangle$ for any 
$v\in{\it l}=T^{\perp}_e(S_{\mathfrak b}\cdot e)$ and $X,Y\in\mathfrak s
=T_e(S_{\mathfrak b}\cdot e)$.  Hence we have 
$A_v\vert_{\mathfrak a\ominus{\it l}}=0$ and 
$A_v\vert_{\mathfrak g_{\lambda}}=\lambda(v){\rm id}$ ($\lambda\in
\triangle_+$), 
where $v\in T^{\perp}_e(S_{\mathfrak b}\cdot e)={\it l}(\subset\mathfrak p)$.  
Therefore, the orbit $S_{\mathfrak b}\cdot e$ is curvature-adapted but it is 
not proper complex equifocal by (ii) of Theorem 1 of [Koi2].  
Hence so are all orbits of the $S_{\mathfrak b}$-action by Lemma 3.3.  

Next we shall show the statement (ii).  
Assume that the root system $\triangle$ of $G/K$ is non-reduced.  
Denote by $A$ the shape tensor of the orbit 
$S_{\mathfrak b,{\it l}_1,\cdots,{\it l}_k}\cdot e\,(\subset\,AN)$.  
Also, let $\xi_0\in\mathfrak b$ and $\xi^i_{t_i}:=
\frac{1}{\cosh(\vert\lambda_i\vert t_i)}\xi^i-\frac{1}{\vert\lambda_i\vert}
\tanh(\vert\lambda_i\vert t_i)H_{\lambda_i}$ ($\xi^i\in
\mathfrak g_{\lambda_i}$) be a unit (tangent) vector of ${\it l}_i$.  Then, 
according to Lemma 3.2, we see that 
$$\begin{array}{l}
\displaystyle{A_{\xi_0}\vert_{\mathfrak s_{\mathfrak b,{\it l}_1,\cdots
{\it l}_k}\cap(\mathfrak a+\sum_{i=1}^k\mathfrak g_{\lambda_i})}=0,}\\
\displaystyle{A_{\xi_0}\vert_{\mathfrak g_{\mu}}=\mu(\xi_0){\rm id}
\quad\,(\mu\in\triangle_+\setminus\mathop{\cup}_{i=1}^k\{\lambda_i\}),}\\
\displaystyle{A_{\xi^i_{t_i}}\vert_{\mathfrak a\ominus
(\mathfrak b+\sum_{j=1}^k{\bf R}H_{\lambda_j})}=0,}\\
\displaystyle{A_{\xi^i_{t_i}}\vert_{{\rm Ker}({\rm ad}(\xi^i)
\vert_{\mathfrak g_{\lambda_i}})\ominus{\bf R}\xi^i}
=-\vert\lambda_i\vert\tanh(\vert\lambda_i\vert t_i){\rm id}}\\
\displaystyle{A_{\xi^i_{t_i}}\vert_{({\bf R}\xi^i+{\bf R}H_{\lambda_i})
\ominus{\it l}_i}=-\vert\lambda_i\vert\tanh(\vert\lambda_i\vert t_i){\rm id}}
\end{array}\leqno{(3.7)}$$
and that, in case of $2\lambda_i\in\triangle_+$, 
$A_{\xi^i_{t_i}}\vert_{{\rm Im}({\rm ad}(\theta\xi^i)
\vert_{\mathfrak g_{2\lambda_i}})+\mathfrak g_{2\lambda_i}}$ has two 
eigenvalues 
$$
\mu_i^+:=-\frac32\vert\lambda_i\vert
\tanh(\vert\lambda_i\vert t_i)
+\frac12\vert\lambda_i\vert\sqrt{2-\tanh^2(\vert\lambda_i\vert t_i)}
$$
and
$$
\mu_i^-:=-\frac32\vert\lambda_i\vert
\tanh(\vert\lambda_i\vert t_i)-
\frac12\vert\lambda_i\vert\sqrt{2-\tanh^2(\vert\lambda_i\vert t_i)}
$$
with the same multiplicity.  Note that $\mathfrak g_{\lambda_i}=
{\rm Ker}({\rm ad}(\xi^i)\vert_{\mathfrak g_{\lambda_i}})\oplus
{\rm Im}({\rm ad}(\theta\xi^i)\vert_{\mathfrak g_{2\lambda_i}})$.  
The eigenspace for $\mu_i^+$ (resp. $\mu_i^-$) is spanned by 
$$\begin{array}{l}
\displaystyle{Z_{\xi^i,Y}^+:=[\theta\xi^i,Y]+\vert\lambda_i\vert
\left(\sinh(\vert\lambda_i\vert t_i)-\sqrt{\sinh^2(\vert\lambda_i\vert t_i)+2}
\right)Y{\rm's}\,\,\,\,(Y\in\mathfrak g_{2\lambda_i})}\\
\displaystyle{{\rm(resp.}\,\,\,\,Z_{\xi^i,Y}^-:=[\theta\xi^i,Y]+
\vert\lambda_i\vert\left(\sinh(\vert\lambda_i\vert t_i)+
\sqrt{\sinh^2(\vert\lambda_i\vert t_i)+2}\right)Y{\rm's}\,\,\,\,
(Y\in\mathfrak g_{2\lambda_i}){\rm))}.}
\end{array}$$
Denote by $R$ the curvature tensor of $\langle\,\,,\,\,\rangle$.  
Also, denote by $X_{\mathfrak f}$ (resp. $X_{\mathfrak p}$) the 
$\mathfrak f$-component (resp. the $\mathfrak p$-component) of 
$X\in\mathfrak g$.  Then we have 
$$\begin{array}{l}
\displaystyle{\left(R(Z^{\pm}_{\xi^i,Y},\xi^i_{t_i})\xi^i_{t_i}
\right)_{\mathfrak p}=-a[[(Z^{\pm}_{\xi^i,Y})_{\mathfrak p},(\xi^i_{t_i})_{\mathfrak p}],
(\xi^i_{t_i})_{\mathfrak p}]}\\
\hspace{3.2truecm}
\displaystyle{=a(-[[Z^{\pm}_{\xi^i,Y},\xi^i_{t_i}],\xi^i_{t_i}]_{\mathfrak p}
+[[(Z^{\pm}_{\xi^i,Y})_{\mathfrak f},(\xi^i_{t_i})_{\mathfrak f}],
(\xi^i_{t_i})_{\mathfrak p}]}\\
\hspace{3.5truecm}
\displaystyle{+[[(Z^{\pm}_{\xi^i,Y})_{\mathfrak f},(\xi^i_{t_i})_{\mathfrak p}],(\xi^i_{t_i})_{\mathfrak f}]
+[[(Z^{\pm}_{\xi^i,Y})_{\mathfrak p},(\xi^i_{t_i})_{\mathfrak f}],
(\xi^i_{t_i})_{\mathfrak f}])}
\end{array}\leqno{(3.8)}$$
for some non-zero constant $a$, where we note that $a=1$ if the metric of 
$G/K$ is induced from the restriction of the Killing form of $\mathfrak g$ 
to $\mathfrak p$.  Also we have 
$$[[(Z^{\pm}_{\xi^i,Y})_{\mathfrak p},(\xi^i_{t_i})_{\mathfrak f}],
(\xi^i_{t_i})_{\mathfrak f}]=0,\leqno{(3.9)}$$
$$[[(Z^{\pm}_{\xi^i,Y})_{\mathfrak f},(\xi^i_{t_i})_{\mathfrak f}],
(\xi^i_{t_i})_{\mathfrak p}]=
-\frac{\tanh(\vert\lambda_i\vert t_i)}{\vert\lambda_i\vert
\cosh(\vert\lambda_i\vert t_i)}
[[[\theta\xi^i,Y]_{\mathfrak f},\xi^i_{\mathfrak f}],H_{\lambda_i}]
\leqno{(3.10)}$$
and 
$$[[(Z^{\pm}_{\xi^i,Y})_{\mathfrak f},(\xi^i_{t_i})_{\mathfrak p}],
(\xi^i_{t_i})_{\mathfrak f}]=
\frac{\vert\lambda_i\vert\tanh(\vert\lambda_i\vert t_i)}
{\cosh(\vert\lambda_i\vert t_i)}[[\theta\xi^i,Y]_{\mathfrak p},
\xi^i_{\mathfrak f}].\leqno{(3.11)}$$
Let $\eta$ (resp. $\bar{\eta}$) be the element of $\mathfrak a+\mathfrak n$ 
with $\eta_{\mathfrak f}=[[\theta\xi^i,Y]_{\mathfrak f},\xi^i_{\mathfrak f}]$ 
(resp. $\bar{\eta}_{\mathfrak p}=[[\theta\xi^i,Y]_{\mathfrak p},
\xi^i_{\mathfrak f}]$).  Then it follows from $(3.8)\sim(3.11)$ that 
$$(R(Z^{\pm}_{\xi^i,Y},\xi^i_{t_i})\xi^i_{t_i})_{\mathfrak p}
=-a[[Z^{\pm}_{\xi^i,Y},\xi^i_{t_i}],\xi^i_{t_i}]_{\mathfrak p}
+\frac{a\vert\lambda_i\vert\tanh(\vert\lambda_i\vert t_i)}
{\cosh(\vert\lambda_i\vert t_i)}(2\eta_{\mathfrak p}+\bar{\eta}_{\mathfrak p}),
$$
that is, 
$$R(Z^{\pm}_{\xi^i,Y},\xi^i_{t_i})\xi^i_{t_i}
=-a[[Z^{\pm}_{\xi^i,Y},\xi^i_{t_i}],\xi^i_{t_i}]
+\frac{a\vert\lambda_i\vert\tanh(\vert\lambda_i\vert t_i)}
{\cosh(\vert\lambda_i\vert t_i)}(2\eta+\bar{\eta}).
\leqno{(3.12)}$$
We have $[\xi^i,\theta\xi^i]=bH_{\lambda_i}$ for some non-zero constant $b$.  
By simple calculation, we have 
$$\begin{array}{l}
\hspace{0.7truecm}
\displaystyle{[[Z^{\pm}_{\xi^i,Y},\xi^i_{t_i}],\xi^i_{t_i}]}\\
\displaystyle{=2\vert\lambda_i\vert\tanh^2(\vert\lambda_i\vert t_i)
\left(-\frac{3b\vert\lambda_i\vert^2}{\sinh(\vert\lambda_i\vert t_i)}
+\sinh(\vert\lambda_i\vert t_i)\mp\sqrt{\sinh^2(\vert\lambda_i\vert t_i)+2}
\right)Y}\\
\hspace{1truecm}
\displaystyle{+\tanh^2(\vert\lambda_i\vert t_i)[\theta\xi^i,Y].}\\
\end{array}
\leqno{(3.13)}$$
From $(3.12)$ and $(3.13)$, it follows that 
$R(Z^{\pm}_{\xi^i,Y},\xi^i_{t_i})\xi^i_{t_i}$ belongs to 
${\rm Im}\,{\rm ad}(\theta\xi^i)\oplus\mathfrak g_{2\lambda_i}$.  Hence 
$R(\cdot,\xi^i_{t_i})\xi^i_{t_i}$ preserves 
${\rm Im}\,{\rm ad}(\theta\xi^i)\oplus\mathfrak g_{2\lambda_i}$ invariantly.  
It is clear that so is also $A_{\xi^i_{t_i}}$.  From $(3.12)$ and 
$(3.13)$, we have 
$[R(\cdot,\xi^i_{t_i})\xi^i_{t_i},A_{\xi^i_{t_i}}]
\vert_{{\rm Im}\,{\rm ad}(\theta\xi^i)\oplus\mathfrak g_{2\lambda_i}}\not=0$, 
under a suitable choice of $t_i$.  
Therefore, 
$S_{\mathfrak b,{\it l}_1,\cdots,{\it l}_k}\cdot e$ is not curvature-adapted 
under suitable choices of ${\it l}_1,\cdots,{\it l}_k$.  
Then, so are all orbits of the 
$S_{\mathfrak b,{\it l}_1,\cdots,{\it l}_k}$-action by Lemma 3.4.  
Furthermore, it follows from Lemma 3.3 that all orbits of the 
$S_{\mathfrak b,{\it l}_1,\cdots,{\it l}_k}$-action are not curvature-adapted 
under arbitrary choices of ${\it l}_1,\cdots,{\it l}_k$.  
Also, it follows from the second relation 
of $(3.7)$ that $S_{\mathfrak b,{\it l}_1,\cdots,{\it l}_k}\cdot e$ 
(hence all orbits of the $S_{\mathfrak b,{\it l}_1,\cdots,{\it l}_k}$-action) 
is not proper complex equifocal in case of $\mathfrak b\not=\{0\}$.  
\hspace{12.1truecm}q.e.d.

\vspace{0.5truecm}

From this proposition, we obtain the statements of Theorem B.  
Also, we have the following fact.  

\vspace{0.5truecm}

\noindent
{\bf Proposition 3.6.} {\sl If $\mathfrak b=\{0\}$, then the 
$S_{\mathfrak b,{\it l}_1,\cdots,{\it l}_k}$-action possesses the only 
minimal orbit.}

\vspace{0.5truecm}

\noindent
{\it Proof.} According to Lemma 3.3, the $S_{\mathfrak b,{\it l}_1,\cdots,
{\it l}_k}$-action is conjugate to $S_{\mathfrak b,\bar{\it l}_1,\cdots,
\bar{\it l}_k}$-action, where $\bar{\it l}_i$ is the orthogonal projection of 
${\it l}_i$ onto $\mathfrak g_{\lambda_i}$.  Hence they are orbit equivalent 
to each other.  Hence we suffice to show that the statement of this 
proposition holds for the $S_{\mathfrak b,\bar{\it l}_1,\cdots,\bar{\it l}_k}$-
action.  Let $\xi^i$ be a unit vector of $\bar{\it l}_i$. Take $p\in AN$.  
We can express as $p=\gamma_{\xi^1}(t_1)\cdots\gamma_{\xi^k}(t_k)$ for some 
$t_1,\cdots,t_k\in{\bf R}$, where $\gamma_{\xi^i}$ is the geodesic with 
$\dot{\gamma}_{\xi^i}(0)=\xi^i$.  Set $\hat{\it l}_i:={\bf R}
\{\frac{1}{\cosh(\vert\lambda_i\vert t_i)}\xi^i
-\frac{1}{\vert\lambda_i\vert}\tanh(\vert\lambda_i\vert t_i)H_{\lambda_i}\}$ 
($i=1,\cdots,k$).  For simplicity, set 
$\xi^i_{t_i}:=\frac{1}{\cosh(\vert\lambda_i\vert t_i)}\xi^i
-\frac{1}{\vert\lambda_i\vert}\tanh(\vert\lambda_i\vert t_i)H_{\lambda_i}$.  
According to the proof of Lemma 3.3, we have 
$$(\gamma_{\xi^1}(t_1)\cdots\gamma_{\xi^k}(t_k))^{-1}
S_{\mathfrak b,\bar{\it l}_1,\cdots,\bar{\it l}_k}(\gamma_{\xi^1}(t_1)\cdots
\gamma_{\xi^k}(t_k))=S_{\mathfrak b,\hat{\it l}_1,\cdots,\hat{\it l}_k}.$$
Hence the orbit $S_{\mathfrak b,\bar{\it l}_1,\cdots,\bar{\it l}_k}\cdot p$ 
is congruent to the orbit 
$S_{\mathfrak b,\hat{\it l}_1,\cdots,\hat{\it l}_k}\cdot e$.  
Denote by $A$ the shape tensor of 
$S_{\mathfrak b,\hat{\it l}_1,\cdots,\hat{\it l}_k}\cdot e$.  
According to Lemma 3.2, we have 
$${\rm Tr}\,A_{\xi^i_{t_i}}=-\vert\lambda_i\vert\tanh(\vert\lambda_i\vert t_i)
\times({\rm dim}\,\mathfrak g_{\lambda\i}+2{\rm dim}\,\mathfrak g_{2\lambda_i})
\quad(i=1,\cdots,k).$$
Hence the orbit $S_{\mathfrak b,\hat{\it l}_1,\cdots,\hat{\it l}_k}\cdot e$ 
is minimal if and only if $t_1=\cdots=t_k=0$, where we note that 
$T^{\perp}_e(S_{\mathfrak b,\hat{\it l}_1,\cdots,\hat{\it l}_k}\cdot e)
={\bf R}\{\xi^1_{t_1},\cdots,\xi^k_{t_k}\}$ because of $\mathfrak b=\{0\}$.  
That is, the orbit $S_{\mathfrak b,\bar{\it l}_1,\cdots,\bar{\it l}_k}\cdot p$ 
is minimal if and only if $p=e$.  Thus the orbit 
$S_{\mathfrak b,\bar{\it l}_1,\cdots,\bar{\it l}_k}$-action posseses the only 
minimal orbit 
$S_{\mathfrak b,\bar{\it l}_1,\cdots,\bar{\it l}_k}\cdot e$.  
This completes the proof.  \hspace{6.7truecm}q.e.d.

\vspace{0.5truecm}

From this proposition, we obtain the statement of Theorem C.  
At the end of this paper, we propose the following question.  

\vspace{0.5truecm}

\noindent
{\bf Question.} {\sl Is any complex hyperpolar action without singular orbit 
on a symmetric space of non-compact type orbit equivalent to 
either the $S_{\mathfrak b}$-action ($\mathfrak b\subset\mathfrak a$) as in 
Example 1 or the $S_{\mathfrak b,{\it l}_1,\cdots,{\it l}_k}$-action 
(${\it l}_i\,:\,$a one dimensional subspace of $\mathfrak g_{\lambda_i}$ 
($i=1,\cdots,k$), $\mathfrak b\subset\mathfrak a\ominus{\rm Span}
\{H_{\lambda_i}\,\vert\,i=1,\cdots,k\}$) as in Example 2 ?}

\vspace{0.8truecm}

\centerline{{\bf References}}

\vspace{0.5truecm}

{\small

\noindent
[B] J. Berndt, Homogeneous hypersurfaces in hyperbolic spaces, Math. Z. 
{\bf 229} (1998) 589-600.

\noindent
[BB] J. Berndt and M. Br$\ddot u$ck, Cohomogeneity one actions on hyperbolic 
spaces, J. Reine Angew. 

Math. {\bf 541} (2001) 209-235.

\noindent
[BT1] J. Berndt and H. Tamaru, Homogeneous codimension one foliations on 
noncompact sym-

metric space, J. Differential Geometry {\bf 63} (2003) 1-40.

\noindent
[BT2] J. Berndt and H. Tamaru, Cohomogeneity one actions on noncompact 
symmetric spaces 

with a totally geodesic singular orbit, 
Tohoku Math. J. {\bf 56} (2004) 163-177.

\noindent
{\small [BV] J. Berndt and L. Vanhecke, 
Curvature adapted submanifolds, 
Nihonkai Math. J. {\bf 3} (1992) 

177-185.

\noindent
[Ch] U. Christ, 
Homogeneity of equifocal submanifolds, J. Differential Geometry 
{\bf 62} (2002) 1-15.

\noindent
[E] H. Ewert, A splitting theorem for equifocal submanifolds in simply 
connected compact symme-

tric spaces, Proc. of Amer. Math. Soc. {\bf 126} (1998) 2443-2452.


\noindent
[G1] L. Geatti, 
Invariant domains in the complexfication of a noncompact Riemannian 
symmetric 

space, J. of Algebra {\bf 251} (2002) 619-685.

\noindent
[G2] L. Geatti, 
Complex extensions of semisimple symmetric spaces, manuscripta math. {\bf 120} 

(2006) 1-25.

\noindent
[HLO] E. Heintze, X. Liu and C. Olmos, 
Isoparametric submanifolds and a 
Chevalley type rest-

riction theorem, Integrable systems, geometry, and topology, 151-190, 
AMS/IP Stud. Adv. 

Math. 36, Amer. Math. Soc., Providence, RI, 2006.

\noindent
[HPTT] E. Heintze, R.S. Palais, C.L. Terng and G. Thorbergsson, 
Hyperpolar actions on symme-

tric spaces, Geometry, topology and physics for Raoul Bott (ed. S. T. Yau), Conf. Proc. 

Lecture Notes Geom. Topology {\bf 4}, 
Internat. Press, Cambridge, MA, 1995 pp214-245.

\noindent
[He] S. Helgason, 
Differential geometry, Lie groups and symmetric spaces, 
Academic Press, New 

York, 1978.

\noindent
[Koi1] N. Koike, 
Submanifold geometries in a symmetric space of non-compact 
type and a pseudo-

Hilbert space, Kyushu J. Math. {\bf 58} (2004) 167-202.

\noindent
[Koi2] N. Koike, 
Complex equifocal submanifolds and infinite dimensional anti-
Kaehlerian isopara-

metric submanifolds, Tokyo J. Math. {\bf 28} (2005) 201-247.

\noindent
[Koi3] N. Koike, 
Actions of Hermann type and proper complex equifocal submanifolds, 
Osaka J. 

Math. {\bf 42} (2005) 599-611.

\noindent
[Koi4] N. Koike, 
A splitting theorem for proper complex equifocal submanifolds, Tohoku Math. 

J. {\bf 58} (2006) 393-417.

\noindent
[Koi5] N. Koike, Complex hyperpolar actions with a totally geodesic orbit, 
Osaka J. Math. {\bf 44} 

(2007) 491-503.

\noindent
[Kol] A. Kollross, A Classification of hyperpolar and cohomogeneity one 
actions, Trans. Amer. 

Math. Soc. {\bf 354} (2001) 571-612.

\noindent
[Ma] A. Malcev, 
On the theory of the Lie groups in the large, Mat. Sb. n. Ser. {\bf 16} 
(1945) 163--

190. (Correction: ibid. {\bf 19} (1946) 523--524).

\noindent
[Mi] J. Milnor, 
Curvatures of left invariant metrics on Lie groups, Adv. Math. {\bf 21} 
(1976) 293--329.

\noindent
[Mo] G.D. Mostow, 
On maximal subgroups of real Lie groups, Ann. Math. {\bf 74} (1961) 
503--517.

\noindent
[PT] R.S. Palais and C.L. Terng, Critical point theory and submanifold 
geometry, Lecture Notes 

in Math. {\bf 1353}, Springer, Berlin, 1988.

\noindent
[S1] R. Sz$\ddot{{{\rm o}}}$ke, Complex structures on tangent 
bundles of Riemannian manifolds, Math. Ann. {\bf 291} 

(1991) 409--428.

\noindent
[S2] R. Sz$\ddot{{{\rm o}}}$ke, Automorphisms of certain Stein 
manifolds, Math. Z. {\bf 219} (1995) 357--385.

\noindent
[S3] R. Sz$\ddot{{{\rm o}}}$ke, Adapted complex structures and 
geometric quantization, Nagoya Math. J. {\bf 154} 

(1999) 171--183.

\noindent
[S4] R. Sz$\ddot{{{\rm o}}}$ke, Involutive structures on the 
tangent bundle of symmetric spaces, 
Math. Ann. {\bf 319} 

(2001), 319--348.

\noindent
[T1] C.L. Terng, 
Isoparametric submanifolds and their Coxeter groups, 
J. Differential Geometry 

{\bf 21} (1985) 79--107.

\noindent
[T2] C.L. Terng, 
Proper Fredholm submanifolds of Hilbert space, 
J. Differential Geometry {\bf 29} 

(1989) 9--47.

\noindent
[T3] C.L. Terng, 
Polar actions on Hilbert space, J. Geom. Anal. {\bf 5} (1995) 129--150.

\noindent
[TT] C.L. Terng and G. Thorbergsson, 
Submanifold geometry in symmetric spaces, 
J. Differential 

Geometry {\bf 42} (1995) 665--718.

\noindent
[W] B. Wu, Isoparametric submanifolds of hyperbolic spaces 
Trans. Amer. Math. Soc. {\bf 331} (1992) 

609--626.
}

\vspace{0.5truecm}

\rightline{Department of Mathematics, Faculty of Science, }
\rightline{Tokyo University of Science}
\rightline{26 Wakamiya Shinjuku-ku,}
\rightline{Tokyo 162-8601, Japan}
\rightline{(e-mail: koike@ma.kagu.tus.ac.jp)}

\end{document}